\documentclass[mathpazo]{cicp}

\usepackage{bm}
\usepackage{xcolor}
\newcommand{\mc}[1]{\mathcal{#1}}
\newcommand{\mb}[1]{\mathbb{#1}}
\newcommand{\vf}{\bm{f}}
\newcommand{\vx}{\bm{x}}
\newcommand{\vy}{\bm{y}}
\newcommand{\vu}{\bm{u}}
\newcommand{\vw}{\bm{w}}
\newcommand{\vW}{\bm{W}}
\newcommand{\vb}{\bm{b}}
\newcommand{\vk}{\bm{k}}
\newcommand{\vB}{\bm{B}}

\newcommand{\vtheta}{\bm{\theta}}
\newcommand{\vX}{\bm{X}}
\newcommand{\vK}{\bm{K}}

\usepackage{float}
\usepackage{subfigure}
\usepackage{algorithm}
\usepackage{algorithmic}
\usepackage{enumitem}

\usepackage{booktabs}
\usepackage{float}
\makeatother

\begin{document}
\title{Semi-Discrete in Time Method for Time-Dependent Equations by Random Neural Basis}


\author{Guihong Wang\affil{1},
      Zheng-An Chen\affil{1}, and Tao Luo\affil{1,2}\comma\corrauth}
\address{\affilnum{1}\ School of Mathematical Sciences,
         Shanghai Jiao Tong University,
         Shanghai 200240, P.R. China \\
          \affilnum{2}\ Institute of Natural Sciences, MOE-LSC, CMA-Shanghai, Shanghai Jiao Tong University, Shanghai 200240, P.R. China}
\emails{{\tt luotao41@sjtu.edu.cn} (T.~Luo), {\tt theodore123@sjtu.edu.cn} (G.~Wang), {\tt zhengan\_chen@sjtu.edu.cn} (Z.~Chen)}

\begin{abstract}
Neural network-based solvers for partial differential equations (PDEs) have attracted considerable attention, yet they often face challenges in accuracy and computational efficiency. In this work, we focus on time-dependent PDEs and observe that coupling space and time in a single network can increase the difficulty of approximation. To address this, we propose a Semi-Discrete in Time Method (SDTM) which leverages classical numerical time integrators and random neural basis (RNB). Additional adaptive operations are introduced to enhance the network's ability to capture features across scales to ensure uniform approximation accuracy for multi-scale PDEs. Numerical experiments demonstrate the framework's effectiveness and confirm the convergence of the temporal integrator as well as the network's approximation performance.
\end{abstract}

\ams{65M50,65Q10,68T07,68U07}
\keywords{partial differential equations, semi-discrete method, random neural basis, adaptive initialization.}

\maketitle

\section{Introduction}
\label{sec:intro}


Partial differential equations (PDEs) play a crucial role in fields such as physics \cite{SommerfeldIntroduction1949PDEiP}, engineering \cite{CiarletFinite2002}, chemistry \cite{SchatzQuantum2002}, and finance \cite{DuffyFinite2006}. Since most PDEs do not admit analytical solutions, their numerical approximation constitutes a central problem in scientific computing. Common approaches, such as finite difference, finite element, and spectral methods, primarily address the spatial discretization. For time-dependent PDEs, while fully discrete formulations are also possible, one of the most widely adopted strategies is the method of lines \cite{KamontNumerical1999HFDIaA}: the spatial derivatives are first replaced by difference expression to obtain a system of ordinary differential equations (ODEs), which is then solved using appropriate time integrators \cite{ErnstHairerSolving1993}. 
These spatial discretization methods can be uniformly interpreted as employing a set of local or global basis functions, where the computation of spatial derivatives essentially corresponds to differentiating the basis functions. A notable feature of such methods is their strong dependence on the underlying mesh discretization. In contrast, neural network–based numerical approaches do not rely on explicit mesh structures; instead, they approximate solutions through parameterized neural networks, which offers potential advantages in handling high-dimensional problems or complex geometries \cite{JingrunChenBridging2022JML}.

The NN-based numerical methods \cite{RaissiPhysicsinformed2019JCP, DongLocal2021CMAME, JingrunChenBridging2022JML} for PDEs have attracted a significant amount of research in the past few years, among which particular mention should be made of Physics-Informed Neural Networks (PINNs) \cite{RaissiPhysicsinformed2019JCP}. PINNs leverage network-based automatic differentiation to train a solution that intrinsically satisfies the governing equations and initial or boundary conditions and are widely recognized for their versatility and ease of implementation, thus providing a general mesh-free framework for solving PDEs. Multiple precision-enhancing strategies have been systematically deployed, including adaptive loss weighting \cite{AnagnostopoulosResidualbased2024CMAME,WangWhen2022JCP}, neural architecture refinement \cite{WangPirateNets2024JMLR}, optimization techniques \cite{WangUnderstanding2021SJSC} and beyond. However, PINNs still face significant criticism concerning accuracy and efficiency limitations compared with classic numerical methods.

Neural network–based approaches typically employ deep neural networks (DNNs) to represent the solutions of PDEs. From this perspective, the network can be regarded as a set of basis functions for approximating the target function. The universal approximation theorem \cite{BarronUniversal1993ITIT} establishes that, as the network width increases, neural networks with analytic activation functions can approximate any smooth target function to arbitrary accuracy. In fact, even for two-layer networks, if the parameters of the hidden layer are fixed—yielding the so-called random feature model—the network still satisfies the universal approximation theorem, albeit with a potentially slower convergence rate \cite{ChenDuality2023}. However, the benefit of fixing the hidden-layer parameters is that it substantially simplifies the training process, as the network optimization reduces to a least-squares problem, for which efficient solvers have been extensively studied. Consequently, a variety of methods based on fixed hidden-layer parameters have been developed, including the Random Feature Model (RFM) \cite{ChenMicromacro2024, JingrunChenBridging2022JML}, Extreme Learning Machine (ELM) \cite{WangExtreme2024CMAME,DongLocal2021CMAME}, and Random Neural Network (RNN) \cite{SunLocal2024JCAM, ShangRandomized2024JEM}. In this work, we regard such methods as a set of network basis functions obtained after random initialization, thus we collectively refer to them as Random Network Basis (RNB) methods. The advantage of these method lies in achieving high accuracy with minimal optimization iterations especially for linear PDEs. A key challenge of the RNB method lies in the fact that the neural network bases are typically global, making it difficult to accurately capture local details in complicate problems. To address this limitation, these methods are often combined with domain decomposition techniques \cite{JingrunChenBridging2022JML, DongMethod2021JCP}, whereby multiple local networks are employed to represent the solution of the equation. For time-dependent PDEs, a commonly adopted strategy is the time marching method \cite{DongLocal2021CMAME, MatteyNovel2022CMAME} in which the temporal domain is partitioned into several consecutive blocks to reduce the complexity of global training.

\begin{figure}[htpb]
    \centering
    \includegraphics[width=0.8\textwidth]{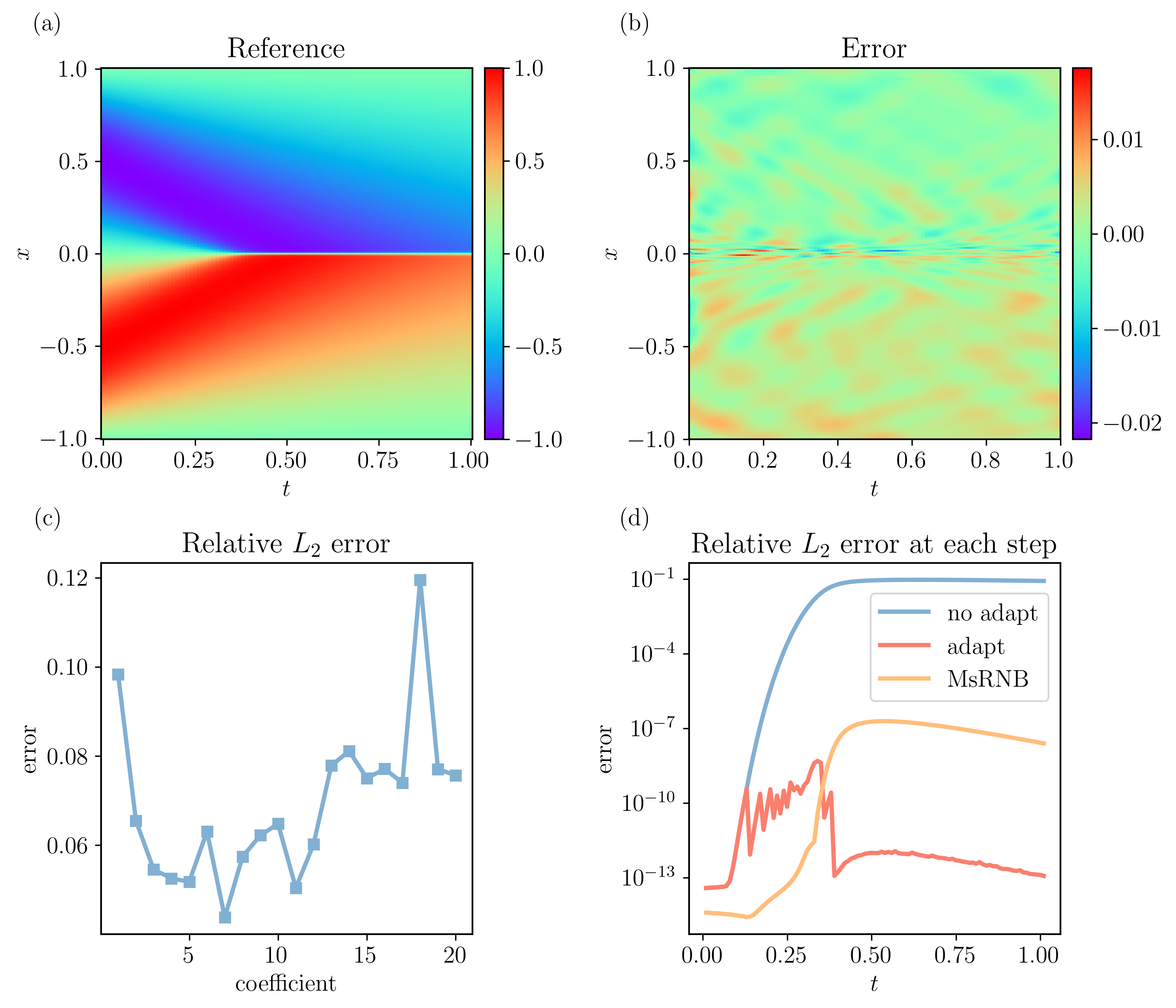}
    \caption{(a) Reference solution for Burgers equation \eqref{eq:bur1d}; (b) Error for DNN of width $[2, 100, 1000, 1]$ in supervised learning after $1\mathrm{e}5$ epochs using Adam; (c) Relative $L_2$ errors in supervised learning for RNB method with different initialization coefficients; (d) Relative $L_2$ errors in supervised learning to each discrete timestep for RNB method of width $[1,1000,1]$ with different initialization method.}
    \label{fig:motivate}
\end{figure}

Going beyond the time marching strategy, an even more thorough approach is to fully discretize the temporal domain, where the neural network is employed solely to represent the solution at each time step as a function of the spatial coordinates. This leads to the so-called optimization-based time integration (OBTI) \cite{ChenTENG2024P4ICML, DonatellaDynamics2023PRA, LuoAutoregressive2022PRL} method, in which the target for next timestep is generated by classic time integrators. Within this framework, the optimization process reduces to approximation problems of the target function at each time step. Common optimizers such as Adam often require multiple iterations, and the additional burden of numerous time discretization steps makes the computation highly time-consuming. To resolve this limitation, time-evolving natural gradient (TENG) method \cite{ChenTENG2024P4ICML} recasts the optimization as a minimal times of least-squares solves. Ref.~\cite{CalabroTime2023AMC} combine the time discretization operation with ELM to solve linear parabolic equations, but the simple addition struggles to solve more complex nonlinear problems. In this paper, we further investigate the optimization process at each time step to preserve high approximation accuracy, yielding a method that remains effective for nonlinear, multi-scale method.

To illustrate our motivation, we conduct a series of experiments from the perspective of accuracy, using reference solutions of the viscous Burgers equation as benchmarks. This equation will be formally introduced in a later section. It should be noted that this benchmark has been widely employed in studies of neural network–based methods \cite{RaissiPhysicsinformed2019JCP, WangPirateNets2024JMLR}. In present experiments, our purpose is not to solve the equation itself, but rather to adopt it in a purely data-driven supervised learning setting, in order to demonstrate the approximation capability of neural networks: even when provided with an exact reference solution, the network still struggles to accurately learn such a challenging mapping. As can be seen in Fig.~\ref{fig:motivate}(a), we show the reference solution, which is computed in traditional spectral method. In Fig.~\ref{fig:motivate}(b), we use a standard DNN to fit the reference with labeled data $\{x_i, u(x_i)\}_{i=1}^{N}$ by Adam. It is evident that the network struggle to accurately or efficiently resolve such solutions under supervised learning paradigms, at least under current setting, let alone within the more demanding physics-informed learning framework. In Fig.~\ref{fig:motivate}(c), we use the same network structure but with fixed hidden-layer parameters of different initialization coefficients and utilize least-squares solvers to train the solution, indicating that it is hard to identify optimal initialization coefficient. The training difficulty in our opinion originates from a fundamental mismatch between the isotropic processing of spatiotemporal coordinates in standard DNN and the anisotropic behavior characteristic of PDE solutions across spatial and temporal dimensions. We then discretize the temporal domain to enable the network to approximate spatial functions exclusively, thereby reducing optimization complexity. In Fig.~\ref{fig:motivate} (d), we employ an ensemble of 101 distinct networks to fit the solution. For the uniform and standard RNB network, the accuracy achieved through supervised learning remains relatively poor for multiscale properties of the equation. However, when an adaptive strategy is incorporated or a multi-scale DNN \cite{LiuMultiscale2020CCP} is employed, the approximation errors remain in a uniform accuracy, outperforming the former two cases. This motivate us to utilize the semi-discrete in time formulation to solve the PDEs, which combines both the advantage of neural networks and the numerical approaches.

The components of our semi-discrete in time method (SDTM) are following:
\begin{enumerate}
    \item \textbf{Random Neural Basis}: Extensive research \cite{JingrunChenBridging2022JML, DongMethod2021JCP,DongLocal2021CMAME, ZhangTransferable2024JSC,ShangRandomized2024JEM, SunLocal2024JCAM} demonstrate that RNBs exhibit exceptional approximation capabilities and scalability across diverse PDE contexts, which preserves the inherent strengths of neural networks: mesh-free computation and automatic differentiation capability. 
    \item \textbf{Time Integration scheme}: This framework enables the implementation of the most prevalent time integration schemes, which is mathematically well-characterized and extensively validated, thereby providing inherited interpretability through documented numerical properties \cite{ErnstHairerSolving1993,StoerIntroduction2002}.    
    \item \textbf{Boundary condition}: Conventional boundary condition enforcement that employs restriction term in loss function introduces iterative errors that compromise stability. We found that the so-called hard boundary constraint strategy \cite{SunSurrogate2020CMAME, DongMethod2021JCP} can mitigate this issue to enhance the accuracy.
    
\noindent\hspace*{-\leftmargini} Moreover, to handle the multi-scale problem, more techniques are required,

    \item \textbf{Multi-scale DNN}: Ref.~\cite{LiMultiscale2020C}  has shown that Multi-scale DNN (MsDNN) outperforms standard DNNs for approximation and solving multiscale elliptic PDEs. Within the RNB method, we identify MsDNN as a specialized initialization strategy that enhances the expressive capacity of the network, enabling it to represent a broader class of target functions.
    
    \item \textbf{Adaptive Initialization Strategy}: Adaptive initialization has long been a key goal for RNB methods, aiming, for instance, to align the spectral properties of the basis functions with those of the target function \cite{JingrunChenBridging2022JML}. In practice, however, only limited information about the solution is available, rendering adaptive strategies largely ineffective. Within the SDTM framework, by contrast, the full information of the target function is accessible, allowing for a straightforward spectral-informed adaptive initialization of the network.
\end{enumerate}
The remainder of this paper are organized as follows: In Section \ref{sec:meth}, we first give a basic introduction to the proposed SDTM framework. To address more challenging problems, Section \ref{sec:adapt} presents strategies for handling multi-scale issues by taking viscous Burgers equation as a motivating example. Section \ref{sec:results} offers several numerical examples to illustrate the effectiveness of the approach and Section \ref{sec:end} concludes the paper with a summary of the results and some further remarks.

\section{Methodology}
\label{sec:meth}
\subsection{Problem Statement}
Consider a time evolving PDE defined on a bounded domain $\Omega \times \mc T \subset \mb R^n \times [0, T]$, let $\vu$ be a function $\Omega \times \mc T \rightarrow \mb R^d$ that satisfies the following initial-boundary value problem
\begin{equation}
\begin{aligned}
  \frac {\partial \vu(\vx,t)}{ \partial t} &= \mc F \vu, \quad (\vx, t) \in \Omega \times \mc T, \\
  \vu(\vx, 0) &= \vu_0(\vx), \\
  \mc B \vu &= 0, \quad \vx \in \partial \Omega,
\end{aligned}
\end{equation}
where $\mc F$ is a linear or nonlinear differential operator, $\mc B$ is a linear boundary operator and $\vu_0$ is the initial condition. We uniformly partition the time domain $[0, T]$ with equidistant time step size $\Delta t$, the OBTI method \cite{ChenTENG2024P4ICML} trains the neural network to approximate the solution of PDE at discrete time steps $t_n = n\Delta t$ as $ \vu_{\vtheta_{t_n}} \approx \vu(\cdot,{t_n}) : \Omega \rightarrow \mb R^d$. The neural network is first optimized to fit the initial condition $\vu_0(\vx)$.  Then network parameters evolve progressively through sequential time steps via directed optimization. 
Take explicit Euler method as an example
\begin{equation}
    \frac{\vu^{n+1} - \vu^{n}}{\Delta t} = \mc F \vu^n.
\end{equation}
Substitute $\vu_{\vtheta_{t_n}}$ into $\vu^n$ above, the target of the next optimization is $\hat \vu_{\vtheta_{t_{n+1}}} =  \vu_{\vtheta_{t_n}} + \Delta t  \mc F \vu_{\vtheta_{t_n}}$. To solve the PDE, the process reduces to $T/\Delta t$ times supervised learning tasks in this case.

\subsection{Random Neural Basis}
A standard $L$ depth fully-connected neural network with $\vx \in \mb R^n$ as input and $\vy \in \mb R^d$ as output can be recursively defined as follows:
\begin{equation}
    \vy_{\vtheta}^{[l]}(\vx) = \vW^{[l]} \cdot \sigma(\vy^{[l-1]}_{\vtheta}(\vx)) + \vb^{l},
\end{equation}
where $\vW^{[l]} \in \mb R^{d_{l}\times d_{l-1}},  \vb^{[l]} \in \mb R^{d_{l}} $, $\sigma$ is the activation function, $d_{l}$ is the width of $l$-th layer, the first layer is $\vy^{[0]}(\vx) = \vx$ and the final layer result is $\vy_{\theta}(\vx) = \vy^{[L+1]}_{\vtheta}(\vx)$. Commonly, $\vtheta$ represents the set of  training parameters $\{\vW^{[L+1]}, \cdots,\vW^{[1]}, \vb^{[L+1]}, \cdots, \vb^{[1]} \}$. Then a supervised task with training data $\{(\vx_i,\hat\vy_i ) \}_{i=1}^{N}$ is formulated as 
\begin{equation}\label{eq:supervised}
    \hat\vtheta =  \arg\min_{\vtheta} \frac{1}{N}\sum_{i=1}^{N} \ell(\vy_{\vtheta}(\vx_i), \hat \vy_i).
\end{equation}
Here, the loss function $\ell(\cdot,\cdot)$ is set as squared $L_2$ norm in this paper, i.e. $l(x,y) = (x-y)^2$. Eq.~\eqref{eq:supervised} is a nonconvex optimization problem thus common optimizers are stochastic gradient-type methods, e.g. SGD, Adam, which typically require multiple iterations to attain a satisfactory model.  Yet if we freeze the hidden-layer parameters, treating the output of network as a linear combination of neural basis, Eq.~\eqref{eq:supervised} becomes a standard linear least-squares problem which can be efficiently and accurately solved by mature numerical solvers.  This approach aligns with RFMs \cite{JingrunChenBridging2022JML,ChenMicromacro2024}, RNNs \cite{ShangRandomized2024JEM,SunLocal2024JCAM} and ELMs \cite{DongLocal2021CMAME,WangExtreme2024CMAME}. In this work, we define the random neural basis (RNB) which spans functional space by neural networks with random initialization and fixed hidden-layer parameters. Then we can express
\begin{equation}
\label{eq:rnb}
  \vy_{\theta} (\vx) = \sum_{i=1}^{d_{L}} \vw^{[L+1]}_i\cdot  y^{[L]}_i(\vx) + \vb^{[L+1]},  
\end{equation}
as a linear combination of RNB, where $\vw^{[L + 1]}_i$ is the $i$-th row vector of $\vW^{[L + 1]}$ and $y^{[L]}_i$ is the $i$-th item of $\vy^{[L]}$. Thereafter $\vtheta$ represents $\{\vW^{L+1}, \vb^{[L+1]} \}$ without causing confusion. 

\subsection{Optimization}
Given the selected collocation points $\vX_{\Omega}$ and $\vX_{\partial \Omega}$, where $\Omega$ is considered as a rectangle domain in the paper. We employ either uniform sampling or Latin Hypercube Sampling (LHS), which is more efficient in higher dimension. After adopting a time integrator and with history solutions $\vu^{j} (j\le n)$ in hand, we can calculate a target solution with the form
\begin{equation}
\vu^{n+1}_{\text{target}}(\vx) = \text{INTG}(\vu^{n-k+1}(x),\cdots,\vu^{n}(x), \vu^{n+1}_{\text{target}}(x)),
\end{equation}
where INTG is an abbreviation for integrator. Since this concept has been extensively studied in numerical analysis \cite{ErnstHairerSolving1993,StoerIntroduction2002}, we provide its detailed formulation in the Appendix. The current solution is needed when an implicit integrator is employed, turning the afterward computational process into an elliptic PDE problem.
The loss function reads
\begin{equation}
    \label{eq:loss}
    L = \frac{1}{|\vX_{\Omega}|}\sum_{\vx \in \vX_{\Omega}} \| \vu_{\theta}(x) - \text{INTG}(\vu^{n-k+1}(x),\cdots, \vu^{n}(x), \vu_{\theta}(x) ) \|_2^2 + \frac{\lambda_\text{bc}}{|\vX_{\partial \Omega}|} \sum_{\vx \in \vX_{\partial \Omega}} \|\mc B \vu_{\theta}(x) \|_2^2. 
\end{equation}
When we use RNB model \eqref{eq:rnb}, the above problem becomes a linear least-squares problem that can be solved in a single step. Moreover, to avoid solving nonlinear problems, we prefer explicit schemes or carefully tuned implicit-explicit (IMEX) \cite{HuangNew2024SJNA} schemes in nonlinear equation case.

\subsection{Boundary condition}
\begin{figure}[htpb]
    \centering
    \includegraphics[width=0.8\textwidth]{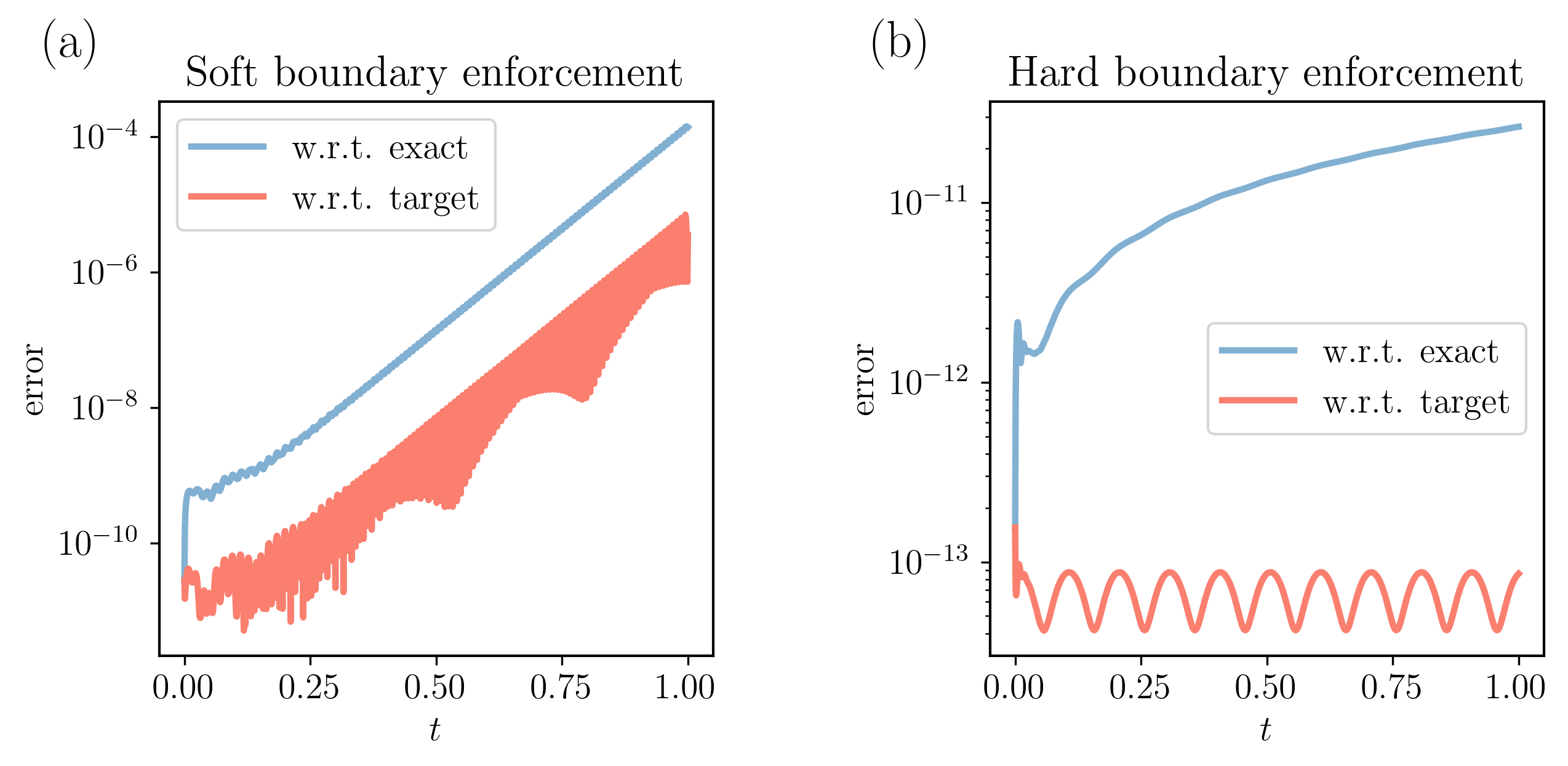}    
    \caption{Different boundary enforcements: We can see from (a) that the soft boundary condition introduces accumulation errors in learning  target function while (b) does not. The blue curves represent the errors with respect to (w.r.t.) exact solutions and the red curves represent the errors w.r.t. the target functions while optimizing the networks.}
    \label{fig:bc}
\end{figure}

In most neural network-based solvers, initial and boundary conditions (BCs) are typically incorporated as additional constraints in the loss function during optimization, as in Eq.~\eqref{eq:loss}, a practice commonly referred to as soft enforcement. However, Ref.~\cite{SunSurrogate2020CMAME} have shown that this enforcement is one of the factors limiting the accuracy of neural network methods. Nevertheless, the SDTM scheme exhibits progressive accumulation of spatial errors during temporal advancement, leading to global accuracy degradation. To mitigate such limitation, Ref.~\cite{DongMethod2021JCP, SunSurrogate2020CMAME} introduces hard boundary enforcement mechanisms designed for Dirichlet and periodic boundary condition to improve solution fidelity, under which we design a solution that satisfies the boundary condition automatically, e.g. periodic Fourier feature layer for periodic boundary condition. The Fourier feature layers read the following form
\begin{equation}
    \Phi(\vx) = \left[ 
    \begin{aligned}
        \sin(\vB \vx) \\
        \cos(\vB \vx)
    \end{aligned}
    \right],
\end{equation}
if we enforce a periodic boundary condition with period $L$, then the items of $\vB$ should be the multiple of $2\pi/L$. For brevity, we omit $2\pi/L$ in the following description of $\vB$. As can be seen in Fig.~\ref{fig:bc}, considering a simple linear advection equation
\begin{equation}
    \label{eq:adv1d}
    u_t + u_x = 0,
\end{equation}
with explicit solution $u(x,t) = \sin(\pi(x - t))$ under periodic boundary condition. When employing RK4 time integrator, the exact solution exhibits high smoothness, enabling ideally accurate approximation of $u_\text{target}$ at each time step. However, soft boundary constraints induce error accumulation through supervised learning cycles in Fig.~\ref{fig:bc}(a), while hard boundary enforcement maintains stable approximation accuracy, confining total error to temporal discretization limits in Fig.~\ref{fig:bc}(b). In our framework, boundary condition treatment is crucial—not only for accuracy but, more importantly, for ensuring overall stability. In this example, soft enforcement progressively degrades the target function, eventually pushing it beyond the network’s representational capacity. Consequently, for the subsequent cases with periodic boundaries, we employ the hard boundary enforcement.

\section{Adaptive Approach for  Multiscale Problem}
\label{sec:adapt}
A primary challenge in RNB method lies in how to initialize the inner layer parameter tailored to heterogeneous target function, which becomes critical in SDTM frameworks for multi-scale problems where solutions evolve high-frequency features during temporal advancement. Consider the 1D viscous Burgers equation \eqref{eq:bur1d} again, the solution becomes sharp during temporal advancement, posing significant challenge to maintain uniform spatial accuracy under fixed RNB. Fig.~\ref{fig:motivate}(d) shows that when the ``shock" takes place at around $t=0.4$, the approximation to $u_\text{target}$ degrades due to emergent high-frequency feature. Thus there are two approaches employed to address this issue. On the one hand, we can find a universal RNB which is suitable for both low-frequency and high-frequency solution. In this paper, we adopt the MsDNN \cite{LiuMultiscale2020CCP, FanMultiscale2019MMS, LiMultiscale2020C} which has good performance for multiscale problems. On the other hand, an adaptive strategy is given that we can reinitialize the network respect to target function with different spectra. To achieve this, we perform a Fourier analysis on the target function and initialize our network according to its frequency content. In this paper, all initializations of the weights are uniform, i.e. $\mc U[-r,r]$, where $r$ denotes the initialization coefficient, and the initializations of bias are standard normal distribution to avoid phase inaccuracy.

\subsection{Multi-scale RNB}
Multi-scale problems represent a foundational challenge in scientific computing, prompting the development of diverse neural architectures \cite{LiuMultiscale2020CCP, FanMultiscale2019MMS} specifically engineered for multi-scale resolution. In this paper we employ MsDNNs which are first proposed in Ref.~\cite{LiuMultiscale2020CCP} and formulated as 
\begin{equation}
    \vy_{\vtheta}^\text{Ms}(\vx) = \vy_{\vtheta} (\vK \odot \vx),
\end{equation}
where $\vK = \{a_1, a_1, \cdots, a_n, a_n\}, a_i \in \mb N^+$, and $|\vK| = d_2$ or $|\vK|=d_3$ for networks with a Fourier feature layer. The MsDNNs help the convergence of high-frequency components of solution in gradient descent optimization according to Ref.~\cite{LiMultiscale2020C,LiuMultiscale2020CCP}. In the case of RNB, we fix the hidden-layer parameters thus the structure of MsDNNs can be seen as a special initialization distribution, thus we refer this variant as multi-scale RNB (MsRNB) in this work. In Fig.~\ref{fig:motivate}(d), we can see that the MsRNB can fit the solution in a uniformly accuracy, outperforming the standard RNB. However, MsRNBs require us to choose a suitable hyperparameter $K$, which is analyzed in Ref.~\cite{HuangFrequencyadaptive2025CMAME} from a frequency perspective. In this paper, we divide $K$ into $n_\text{max}$ equal segments and assign values from 1 to $n_\text{max}$ sequentially. However, we do not investigate in depth how this hyperparameter should be set but which adapts with slight reference to the subsequent approach.

\subsection{Adaptive Initialization Strategy}
\begin{figure}[htpb]
    \centering
    \includegraphics[width=\textwidth]{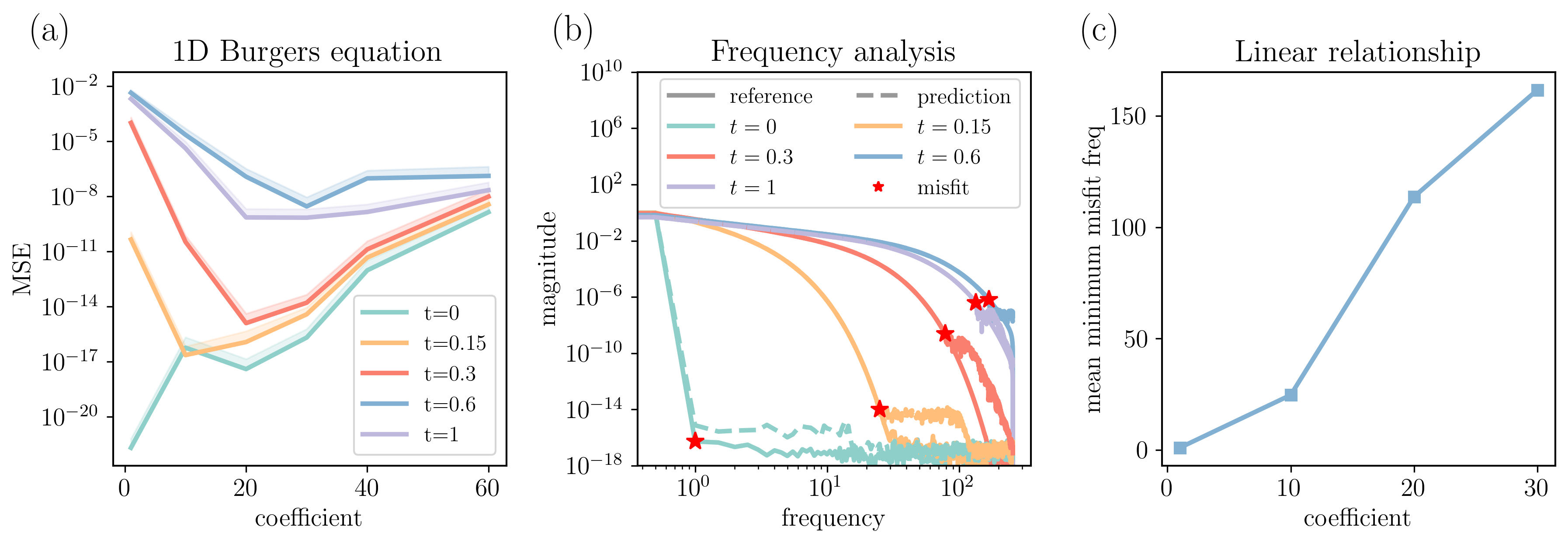}  
    \caption{(a) Various initialization coefficients performance for different target functions; (b) Frequency analysis of the optimal coefficients in (a); (c) Relationship between max misfit frequencies and initialization coefficients.}
    \label{fig:adapt}
\end{figure}

The concept of preconditioning is common in scientific computing, primarily involving the adaptive adjustment of algorithms based on prior information. In RNB methods, Ref.~\cite{ZhangTransferable2024JSC} also employs auxiliary functions to tune the optimal hyper-parameter of the network, aiming to obtain a general RNB which is suitable for a large range of target functions. In this work, we propose a adaptive strategy that can initialize the network according to target functions based on their properties. More precisely, we reinitialize the RNB respect to the frequencies of target functions while the approximation under-performs. As can be seen in Fig.~\ref{fig:motivate}(d), the approximation of supervised learning can destroy the overall performance, indicating that the standard RNB is not suitable. To investigate the performance of RNB under different initialization, we plot the Mean Squared Error (MSE) for varying $r$ in Fig.~\ref{fig:adapt}(a), respect to different target functions. These functions are reference solutions at different time steps of Eq.~\eqref{eq:bur1d}. We use a three-layer network with structure $[1,2,1000,1]$, including a Fourier feature layer of $\vB=[1]$. The results are obtained by averaging over $10$ independent runs. The shaded area indicates the positive standard deviation. Since the plots are drawn on a logarithmic scale, the negative deviation is not displayed in order to avoid undefined values. The optimal coefficient settings vary across different target functions, which suggests that we need to enhance the initialization in order to better fit the subsequent solutions, at least for example \eqref{eq:bur1d}. To explore why the network struggles to learn the target function, a frequency analysis in Fig.~\ref{fig:adapt}(b) reveals that it cannot accurately capture high-frequency components, owing to its inherent limitations in representing such frequencies. The RNB in Fig.~\ref{fig:adapt}(b) is configured with the optimal coefficients as indicated in Fig.~\ref{fig:adapt}(a).  Focused on the minimum misfit frequencies, which are marked as red stars in Fig.~\ref{fig:adapt}(b), the magnitudes of these frequencies are consistent with the precision shown in the left figure, according the the Parseval`s identity. This implies that the magnitudes of minimum misfit frequencies have a direct effect on the approximation errors. In Fig.~\ref{fig:adapt}(c), we plot the relationship between the minimum misfit frequencies and initialization coefficients, the results are averaged over $10$ independent runs, too. We can discover that there is a quasi-linear relationship in Fig.~\ref{fig:adapt}(c), indicating that we can initialize the network according to frequencies of target functions. 

\noindent\textbf{Adaptive Initialization Strategy based on Frequency:}
Taking the one-dimensional case for illustration, for the Fourier transform $\hat{u}(\omega)$ of the target function, we at least expect the function to be learned up to an $L_2$ accuracy of $\epsilon$. Therefore, any frequency component with magnitude greater than $\epsilon$ should be included in the RNB, and based on this criterion we initialize the network. Here, we choose
\begin{equation}
    r = \arg\max_{\omega} | |\hat u(\omega)| \ge \epsilon|,
\end{equation}
then the initialization coefficients should be $\min\{r, r_\text{max}\}$, where the $r_{\max}$ is a hyper-parameter to avoid unstable evolution. In the case of employing a Fourier feature layer, the rescaling $r = r / \max(\vB)$ is required.

To further illustrate this phenomenon, we analyze the frequency characteristics of the basis functions $f_{\bm{k}}(\bm{x}) = \tanh(\sum_{i=1}^n k_i \sin(x_i))$. 
Specifically, we characterize their frequency support, which we define for a given tolerance $\varepsilon$ as the smallest indices $\bm{m}$ beyond which all Fourier coefficients are bounded by $\varepsilon$:
\begin{equation}
    S_{k} = \inf \left\{ \|\bm{m}\|_1 \in \mathbb{N} \ | \  |\hat{f}_{\bm{k}} (\bm{j})| \le \varepsilon,\ \ \forall \| \bm{j}\|_1 >  \|\bm{m}\|_1 \right\}.
\end{equation}
As the parameter $k$---which corresponds to the initialization scale---increases, the frequency support $S_k$ grows. 
The following proposition establishes that this growth is at most linear. See Appendix for the detailed proof.
\begin{proposition}
    Let $f_k(x) = \tanh(\sum_{i=1}^n k_i \sin(x_i))$. The $\varepsilon$-truncated frequency support $S_k$ is bounded by a linear function of the initialization scale $k$:
    \begin{equation}
        S_k \le C \| \bm{k}\|_1,
    \end{equation}
    where the constant $C$ depends on $\varepsilon$ but is independent of $k$.
\end{proposition}

Finally, we present the SDTM algorithm with adaptive initialization in Algorithm 1. For SDTM without reinitialization at each timestep i.e. using standard RNB or MSDNN network, the procedure is similar and trivial to Algorithm 1. We have to comment that the target function \eqref{eq:target} can be computed explicitly while using an explicit integrator. For an implicit integrator, we can calculate $u^n_\text{target}$ in practice that includes trainable parameters using the python package \textit{torch}. Then we subscribe $u^n_\text{target}$ to Eq.~\eqref{eq:lstsq} to generate a least-squares system. In this paper, we mainly focus on explicit time integrators and some results for implicit schemes are presented in Appendix. 
\begin{algorithm}[h]
    \label{alg:adapt}
    \caption{SDTM with adaptive initialization}
    \begin{algorithmic}
        \REQUIRE The initial condition $\vu_i (i < k)$, collocation points $\vX_{\Omega}$, $\vX_{\partial\Omega}$, $k$-step time integrator INTG, RNB model $\vu_{\vtheta}$ with out-layer parameters $\{\vtheta_j\}_{j=1}^{M} := \vtheta$, time steps $t_n (n \le N)$, initialization coefficient $r, r_\text{max}$, tolerance parameter $\epsilon$.
    
        \ENSURE The prediction $\vu^{n}_{\vtheta} (\vx)$ for arbitrary points $\vx\in \vX_\text{test}$.
        
\noindent\textbf{Initialization:} initialize $\vu_{\vtheta}$ with $\mc U[-r, r]$, calculate the predictions $\vu_{\vtheta}^{i} (i<k)$ by solving
$$\arg\min_{\vtheta} \sum_{\vx\in \vX_{\Omega}} \|\vu_{\vtheta}^i(\vx) - \vu_i(\vx) \|_2^2$$
\FOR{$n = k:N$}
\item compute the target function
\begin{equation}
\label{eq:target}
\vu^n_\text{target} (\vx) = \text{INTG}(\vu^{n-k}_{\vtheta}(\vx),\cdots,\vu^{n-1}_{\vtheta}(\vx),\vu^{n}_\text{target}(\vx) ),    
\end{equation}

\IF{$\|\vu^{n-1}_{\vtheta} -\vu^{n-1}_\text{target} \|_2 > \epsilon$}
\item perform the Fourier transform towards $\vu^n_\text{target}$ as $\hat \vu(\vw) = \mc F(\vu^n_\text{target})$, then initialize the hidden-layer parameters of $\vu_{\vtheta}^{n}$ with $\mc{U}[-r, r]$, where
\begin{equation}
    \label{eq:adapt}
    r = \min\{\arg\max_{|\vw|} | |\hat \vu(\vw)| \ge \epsilon|, r_\text{max}\}.
\end{equation}
\ENDIF
\item solve the least-squares problem
\begin{equation}
    \label{eq:lstsq}
    \arg\min_{\vtheta} \frac{1}{|\vX_{\Omega}|}\sum_{\vx\in\vX_{\Omega}}\|\vu^{n}_{\vtheta}(\vx) - \vu^n_\text{target} (x)  \|_2^2 + \frac{\lambda_{\text{bc}}}{|\vX_{\partial \Omega}|}\sum_{\vx\in\vX_{\partial\Omega}}\|\mc B \vu^n_{\vtheta}(\vx)  \|_2^2.    
\end{equation}

\ENDFOR
    \end{algorithmic}
\end{algorithm}

\section{Numerical Results}
\label{sec:results}
In this section, several numerical problems are presented to show the effectiveness and accuracy of SDTM. As higher-order derivatives are obtained via automatic differentiation, we incorporate an $L_2$ penalty with weight $\lambda$ into the linear least-squares solution to reduce the effect of overfitting on derivative accuracy, with $\lambda = 1\mathrm{e}\text{-20}$ in all experiments. For the examples without analytic solution, the reference is solved using an FFT-based pseudospectral method  with a RK4 integrating factor time-stepping scheme. To evaluate the performance, we use relative $L_2$ error and $L$-infinity error as criteria 
\begin{equation}
    e^\text{rel}(u, u^{*}) = \frac{\|u - u^{*}\|_2}{\|u^*\|_2}, \quad e^\text{abs}(u,u^{*}) = \|u-u^{*}\|_\infty.
\end{equation}
\subsection{1D Advection equation}
\begin{figure}[htpb]
    \centering
    \includegraphics[width=\textwidth]{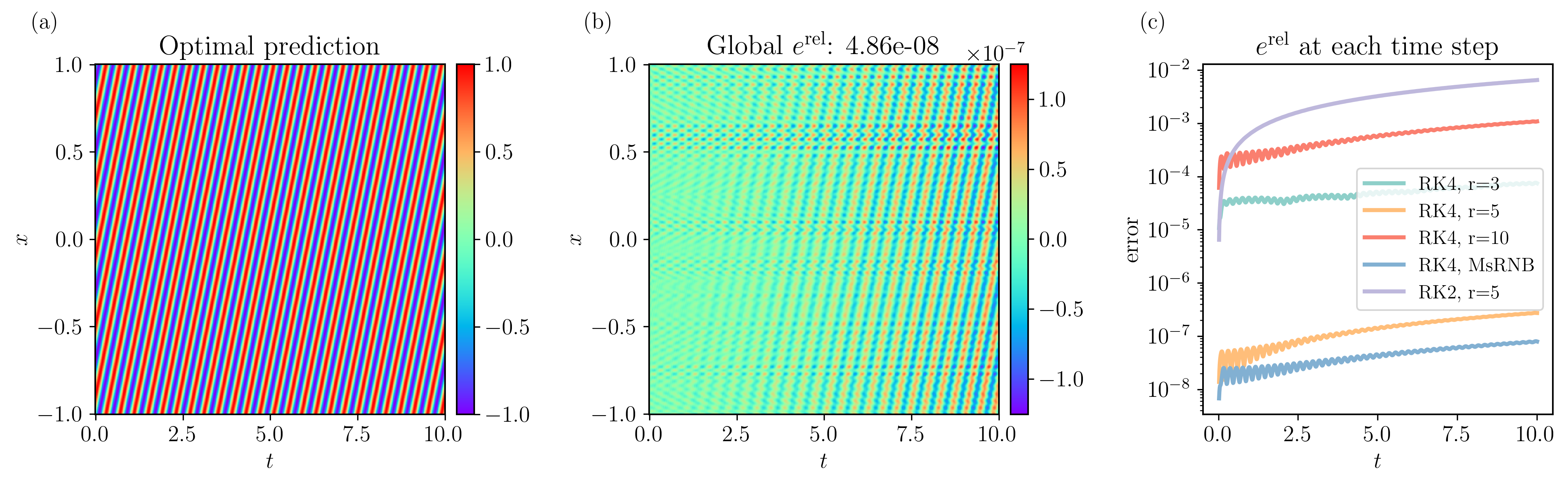}  
    \caption{Advection Equation: (a) Prediction of MsDNN; (b) Error of the prediction; (c) Relative $L_2$ error in each timestep for different networks with RK2 or RK4 integrator.}
    \label{fig:adv1d}
\end{figure}

We begin with a simple linear problem, which prior results have shown to exhibit poor performance in the absence of additional techniques. Consider the linear advection equation \eqref{eq:adv1d}, with exact solution $u(x) = \sin(5\pi (x - t))$, which indicates us that the optimal initialization coefficient is $5$. Next, we employ a three-layer network with widths of $[1,2,100,1]$, which includes one periodic Fourier layer of $\vB=[1]$, to compute the solution of this equation up to $t=10$. We use LHS method to sample $1000$ random collocation points for training and test in a $513 \times 1001$ grid. We use different initialization coefficients of $3, 5, 10$ using RK4 integrator with $\Delta t = 0.001$ and plot the relative $L_2$ error to the exact solution in Fig.~\ref{fig:adv1d}(c), from which we can find that the best performance among the three coefficients is $r=5$, aligning with our adaptive criteria. Moreover, we employ a MsDNN with $n_\text{max}=10$ in the same setting, outperforming the former networks in this case. This optimal prediction is plotted in Fig.~\ref{fig:adv1d}(a), and the error is plotted in Fig.~\ref{fig:ac1d}(b), achieving the global relative $L_2$ error to $4.86\mathrm{e}\text{-8}$. This accuracy is attributed to the suitability of the network as well as the high-precision time integrator. When we instead employ the RK2 integrator, which is again shown in Fig.~\ref{fig:adv1d}(c), we observe that the accuracy is degraded.

\subsection{1D viscous Burgers equation}
\begin{figure}[htpb]
    \centering
    \includegraphics[width=\textwidth]{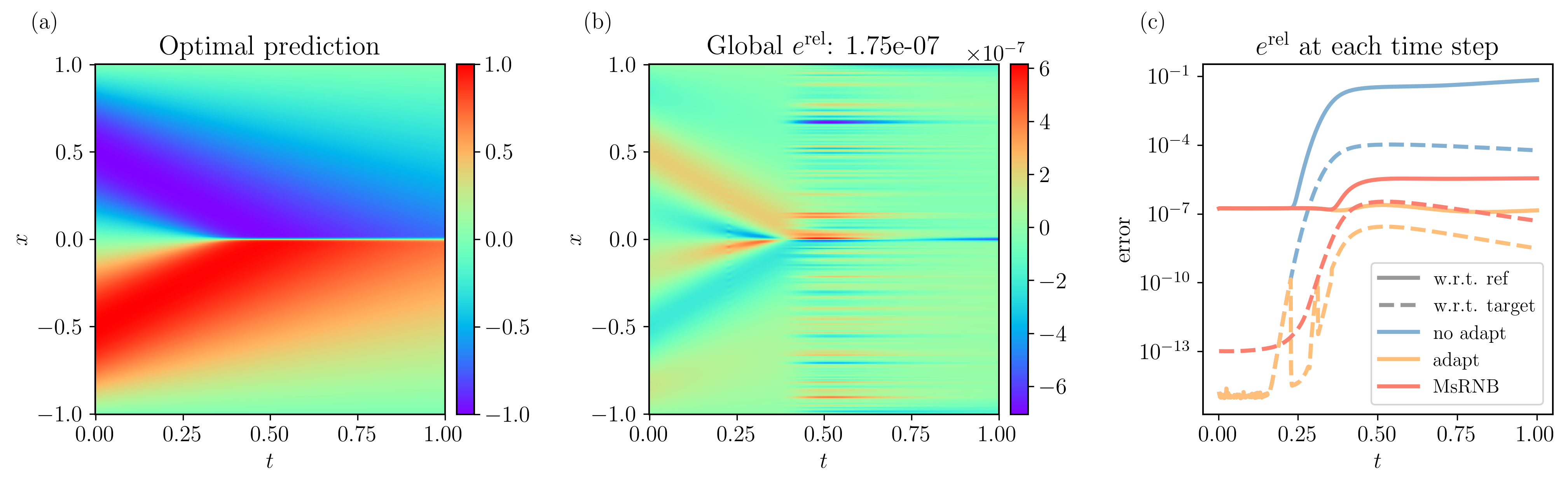}  
    \caption{Burgers Equation: (a) Prediction of adaptive initialized RNB; (b) Error of the prediction; (c) Relative $L_2$ error at each timestep for three different settings. Blue: normal constant initialized RNB, yellow: adaptive initialized RNB, red: MsRNB. The solid curves represent the errors with respect to reference solutions and the dash curves represent the errors to the target functions generated at each timestep.}
    \label{fig:bur1d}
\end{figure}

Consider the 1D viscous Burgers equation,
\begin{equation}
    \label{eq:bur1d}
    \begin{aligned}
    u_t + uu_x &= \nu u_{xx},\quad  (x, t) \in [-1, 1] \times [0, 1] \\
    u(x,0) &= u_0(x),
\end{aligned}
\end{equation}
with periodic boundary condition, where $\nu = 0.01/\pi$, $u_0(x) = -\sin(\pi x)$. The viscous Burgers equation is a fundamental nonlinear PDE that combines advection, diffusion, and shock formation phenomena. It serves as a simplified model for fluid dynamics, turbulence, and wave propagation. This example is popular for PINN methods \cite{WangPirateNets2024JMLR, RaissiPhysicsinformed2019JCP}, showing the difficulty to learn the ``shock" solution accurately. Thus this issue recurs throughout the motivation and analysis of this paper, where we show that a MsRNB or an adaptive initialized RNB can capture solutions for each timestep. In this problem, we use BDF-type integrator with $\Delta t = 1\text{e}\text{-4}$ which is proposed in Ref.~\cite{HuangNew2024SJNA} and select an explicit version which is discretized at $t^{n+2}$. More details are available in Appendix. The BDF-type methods are more efficient than higher-order one-step RK-type methods, which is illustrate in Sec.\ref{sec:cost}. To have better performance, we use a uniform grid of size $4097$ for $x\in[-1,1]$ and test in a $513 \times 1001$ grid in spatiotemporal domain. We show the experimental results in Fig.~\ref{fig:bur1d}, what we are more concerned with are the approximate errors of the networks with respect to the target functions at each step. We test three different settings, the normal constant initialization with $r=1$, the adaptive initialized RNB and MsRNB. We use a three-layer network with widths of  $[1,4,1000,1]$ including a periodic Fourier feature layer of $\vB=[1,2]$. In Fig.~\ref{fig:bur1d}(c), we plot the relative $L_2$ errors to reference solutions with solid curves and plot errors to target functions at each step with dash curves in consistent colors. We can find that the blue and yellow dash curves coincide before $t=0.2$ since we use the same random seed and the networks are the same during this period. The approximate error keeps increasing, thereby deteriorating the overall accuracy of blue solid curve. While in yellow dash curve, we reinitialize the network adaptively using Eq.~\eqref{eq:adapt}, thus the approximate error keeps in a uniform accuracy, contributing to a steady overall error accumulation in yellow solid curve. As for MsRNB, although it also experiences an increase in approximate error during ``shock" formation, it is able to maintain the error at a relatively low level, thereby yielding a satisfactory result. Similar to the former example, we plot the optimal prediction and error in Fig.~\ref{fig:bur1d}(a)(b), where we can see that the global relative $L_2$ error of adaptive initialization reaches $1.75\mathrm{e}\text{-7}$.

\subsection{1D Allen--Cahn equation}
\begin{figure}[htpb]
    \centering
    \includegraphics[width=0.8\textwidth]{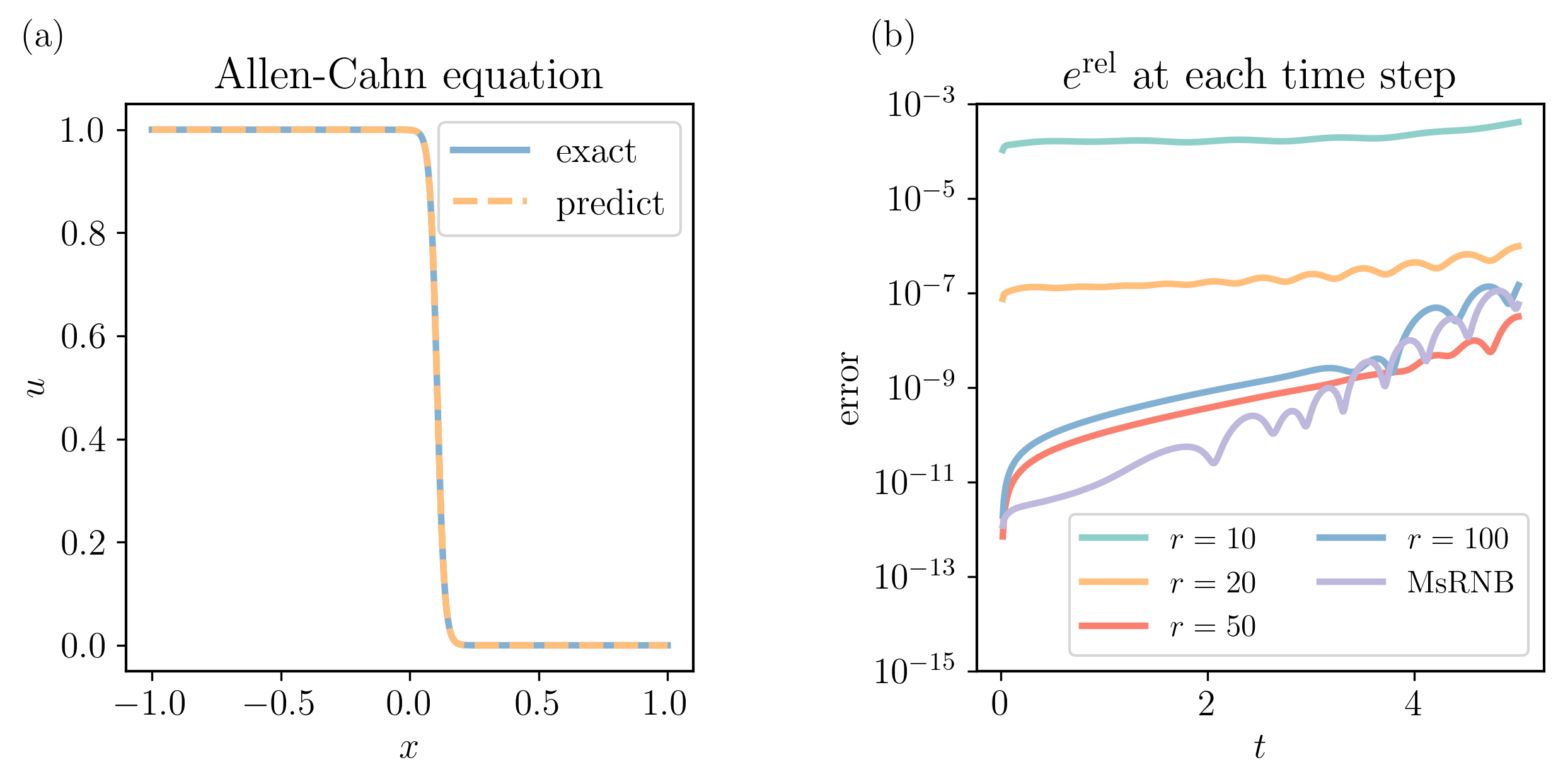}  
    \caption{Allen--Cahn Equation \eqref{eq:ac1d2}: (a) Prediction at $t=5$ of the case $r=50$; (b) Relative $L_2$ error at each timestep for 5 different setting, $r = 10, 20, 50, 100$ and MsRNB.}
    \label{fig:ac1d2}
\end{figure}

The Allen--Cahn equation is a cornerstone of the phase-field modeling methodology, which is used to describe the evolution of micro-structures and interfacial dynamics. We begin with the well-known traveling wave solution \cite{PoochinapanNumerical2022AMC} for 1D Allen--Cahn equation
\begin{equation}
    \label{eq:ac1d2}
    \begin{aligned}
    u_t = \epsilon^2 u_{xx} - u^3 + u,\\
    u(x,0) = u_0(x),
    \end{aligned}
    \end{equation}
the solution has the explicit expression
\begin{equation}
    u(x,t) = \frac{1}{2}\left(1-\tanh(\frac{x-st}{2\sqrt{2}\epsilon})\right), 
\end{equation}
where $s=3\epsilon/\sqrt{2}$ is the speed of the traveling wave. The initial condition is given by the exact solution and in this case the boundary conditions are 
\begin{equation}
    u(-1, t) = 1, \quad u(1,t) = 0.
\end{equation}
In this example, the equation is tested for $\epsilon=0.01$ over $[-1, 1]\times[0, 5]$ using a two-layer network of width $1000$ and $\Delta t = 0.001$. The time integrator and the dataset setting are the same as that of former Burgers equation. Since the target function in this example involves large frequency values across the spectrum, we prefer larger initialization coefficients. We test $r = 10, 20, 50, 100$ in Fig.~\ref{fig:ac1d2}(b), showing that the cases of $r = 50$ and $r=100$ outperform, between which $r =50$ yields a smoother error curve. We attribute this effect to larger initialization, which may introduce higher-frequency perturbations into the derivatives obtained by automatic differentiation. Moreover, we test a MsRNB case with $n_\text{max}=100$, which yields comparable performance. In Fig.~\ref{fig:ac1d2}(a), we plot the prediction at $t=5$ with respect to exact solution. This figure shows that we can fit this sharp solution well and the global relative $L_2$ error is $5.73\mathrm{e}\text{-9}$.  

\begin{figure}[htpb]
    \centering
    \includegraphics[width=\textwidth]{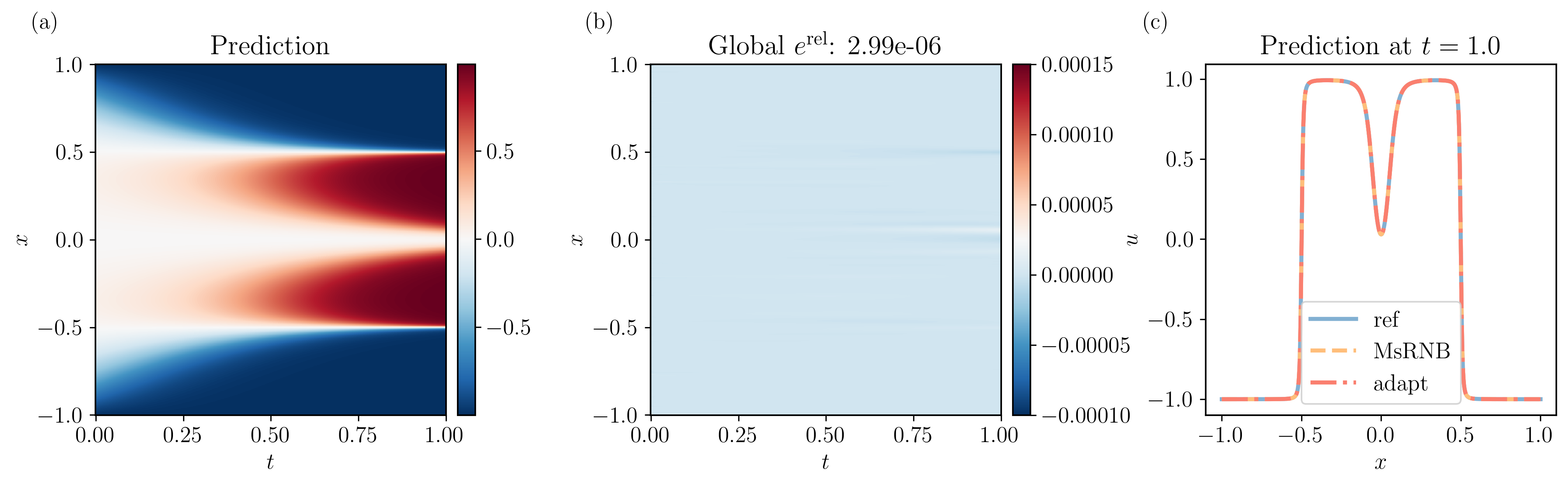}  
    \caption{Allen--Cahn Equation \eqref{eq:ac1d}: (a) Prediction of MsRNB; (b) Error of the prediction; (c) Prediction at $t=1$ under the settings of MsRNB and adaptive initialized RNB.}
    \label{fig:ac1d}
\end{figure}

We consider another 1D Allen--Cahn example which is transformed to the form
\begin{equation}
    \label{eq:ac1d}
    u_t = \epsilon^2 u_{xx} - 5u^3 + 5u, \quad (x,t)\in [-1,1]\times[0,1],
\end{equation}
with periodic boundary condition. The initial condition is $u_0(x) = x^2\cos(\pi x)$ and $\epsilon = 0.01$. This example is also popular for PINN type methods \cite{WangPirateNets2024JMLR, RaissiPhysicsinformed2019JCP}. We use a three-layer network with widths of $[1,4,1000,1]$ including a periodic Fourier feature layer of $\vB=[1,2]$ and employ the same explicit BDF4 method with $\Delta t = 1 \mathrm{e}\text{-4}$. We use LHS method to sample $4000$ random collocation points and test in a $513\times1001$ gird. The adaptive initialized RNB and MsRNB are tested separately which yield comparable performance. In Fig.~\ref{fig:ac1d}(a)(b), we plot the prediction and the error to reference solution of MsRNB setting, which has a global relative $L_2$ error of $2.99\mathrm{e}\text{-6}$. We also plot the predictions of the solution at $t=1$ for the two settings, showing that both can capture the multiscale structure of solution.  

\subsection{2D Allen--Cahn equation}
\begin{figure}[htpb]
    \centering
    \includegraphics[width=\textwidth]{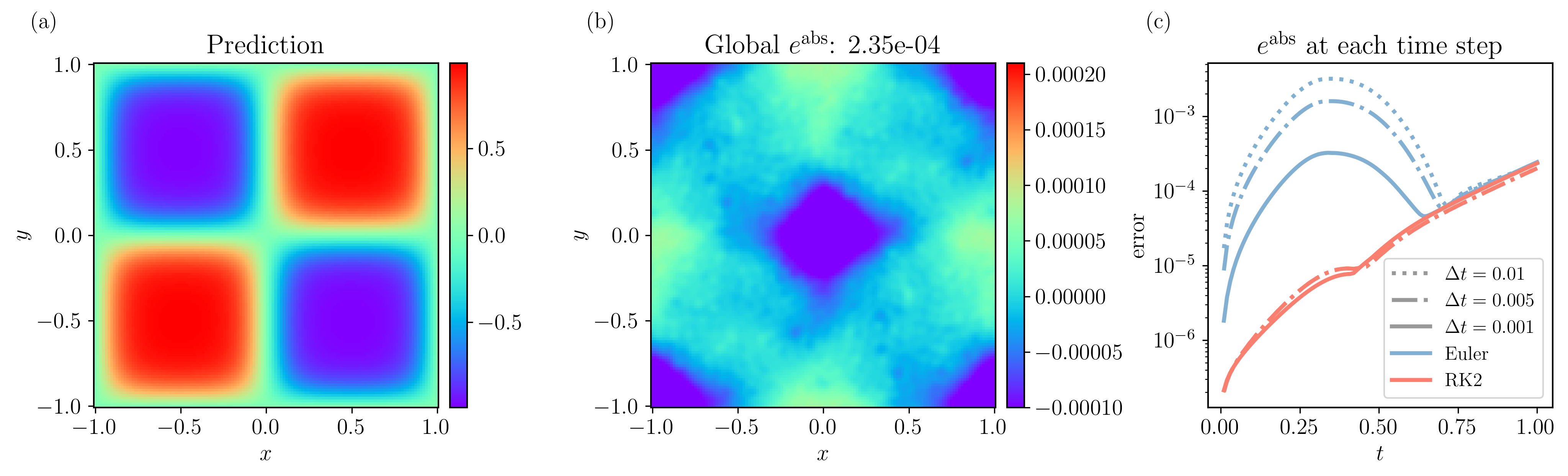}  
    \label{fig:ac2d}
    \caption{2D Allen--Cahn Equation \eqref{eq:ac2d}: (a) Prediction of RK2 with $\Delta t = 0.001$; (b) Error of the prediction; (c) L-infinity errors to each timestep for Euler and RK2 integrator with $\Delta t = 0.01, 0.005, 0.001$.}
\end{figure}

In the latter two examples, we focus on 2D cases to show effectiveness of SDTM for higher dimensions. Consider the 2D Allen--Cahn equation which have the similar expression to the 1D case
\begin{equation}
    \label{eq:ac2d}
    u_t = \epsilon^2 \Delta u - u^3 + u, \quad (x,y,t)\in [-1,1]^2\times [0,10], 
\end{equation}
the initial condition is $u_0(x,y) = 0.05\sin(\pi x)\sin(\pi y)$, $\epsilon=0.1$, and we use a periodic boundary condition. We use LHS method to sample $4000$ random collocation points for training and test in a $128\times128\times101$ grid. We plot the solution at $t=10$ in Fig.~\ref{fig:ac2d}(a), we observe that the final steady forms four distinct phase domains. In fact, this example reaches equilibrium at $t=4$ and stays steady afterwards. We test for Euler and RK2 integrator for $\Delta t = 0.01, 0.005, 0.001$ and use RNB with widths of $[2, 8, 1000,1]$ containing a Fourier feature layer of $B=[1,2,3,4]$. We can see from Fig.~\ref{fig:ac2d}(c) that all the settings converge to the same error. The performance of RK2 and $\Delta t=0.001$ is plotted in Fig.~\ref{fig:ac2d}(a)(b), where we can see that the $L$-infinity error is $2.35\mathrm{e}\text{-4}$.

\subsection{2D Navier--Stokes equation}
\begin{figure}[htpb]
    \centering
    \includegraphics[width=0.8\textwidth]{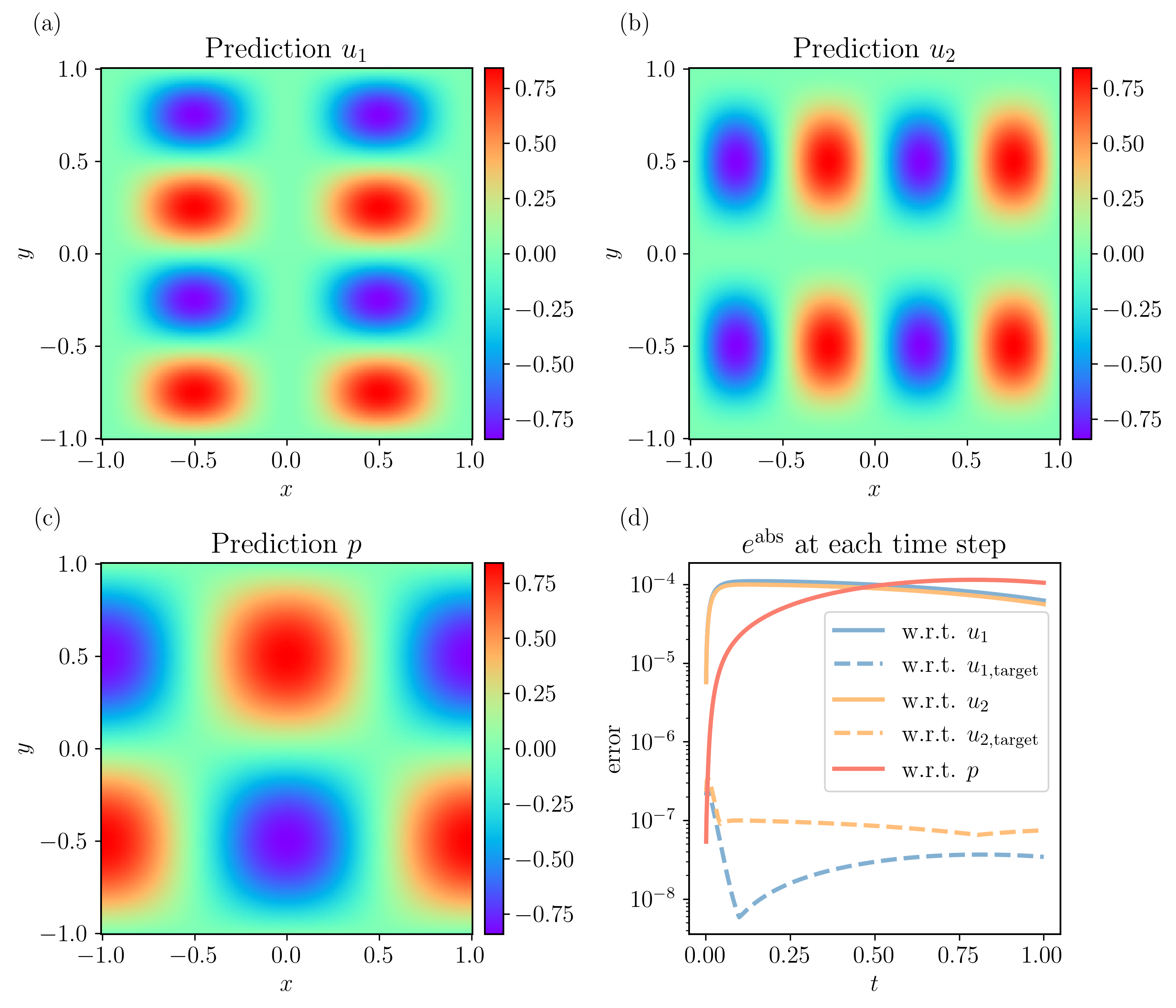}  
    \label{fig:ns2d}
    \caption{Navier--Stokes equation: (a)(b)(c): Predictions for $u_1$,$u_2$ and $p$; (d):$L$-infinity errors at each timestep for $u_1, u_2$ and $p$.}
\end{figure}

Finally, we consider 2D Navier--Stokes equation which reads
\begin{equation}
    \label{eq:ns2d}
    \begin{aligned}
        \frac{\partial \vu}{\partial t} + \vu \cdot\nabla\vu &-\nu\Delta\vu +\nabla p = \vf, \\
        \nabla \cdot \vu &= 0,
    \end{aligned}
\end{equation}
where $\vu$ is velocity and $p$ is the pressure. We choose $\nu = 1$ and the exact solution are given by Ref.~\cite{HuangStability2023SJNAa} as 
\begin{equation}
    \label{eq:ns2d-exact}
    \begin{aligned}
        u_1(x,y,t) &= \sin(2\pi y)\sin^2(\pi x) \sin(t),\\
        u_2(x,y,t) &= \sin(2\pi x)\sin^2(\pi y) \sin(t),\\
        p(x,y,t) &= \cos(\pi x)\sin(\pi y)\sin(y).
    \end{aligned}
\end{equation}
The initial condition is given by Eq.~\eqref{eq:ns2d-exact} using no-slip boundary condition that $\vu=0$ on $\partial \Omega$. To solve this problem, we perform divergence operations on the first equation of Eq.~\eqref{eq:ns2d}. Since $\vu$ is divergence free, we can get
\begin{equation}
    \label{eq:ns2d-pressure}
    \Delta p = -\left(\frac{\partial u_1} {\partial x}\frac{\partial u_1} {\partial x}+ 2\frac{\partial u_1} {\partial y}\frac{\partial u_2} {\partial x}+\frac{\partial u_2} {\partial y}\frac{\partial u_2} {\partial y}\right) - \nabla\cdot \vf,
\end{equation}
which is a Poisson equation for pressure $p$. The calculating procedure is that we first compute $\vu^{n+1}$ using $\vu^{n}$, $p^{n}$ with SDTM and then solve Eq.~\eqref{eq:ns2d-pressure} with RNB method, which can be solved in a single iteration. The computation cost at each timestep is to solve two linear least-squares problems. We use a RNB with widths of $[2, 8, 1000, 2]$ to express $\vu$ which include a Fourier layer, and use another network with widths of $[2,8,1000,1]$ to represent $p$. We also use the explicit BDF4 integrator with $\Delta t = 1\mathrm{e}\text{-4}$. The predictions are plotted in Fig.~\ref{fig:ns2d}(a)(b)(c), and the $L$-infinity errors at each timestep are plotted in Fig.~\ref{fig:ns2d}(d). The global relative $L_2$ errors of $u_1$, $u_2$, $p$ are respectively $1.64\mathrm{e}\text{-4}$, $1.61\mathrm{e}\text{-4}$, $1.05\mathrm{e}\text{-4}$.  

\subsection{Convergence analysis}
\begin{figure}[htpb]
    \centering
    \includegraphics[width=0.8\textwidth]{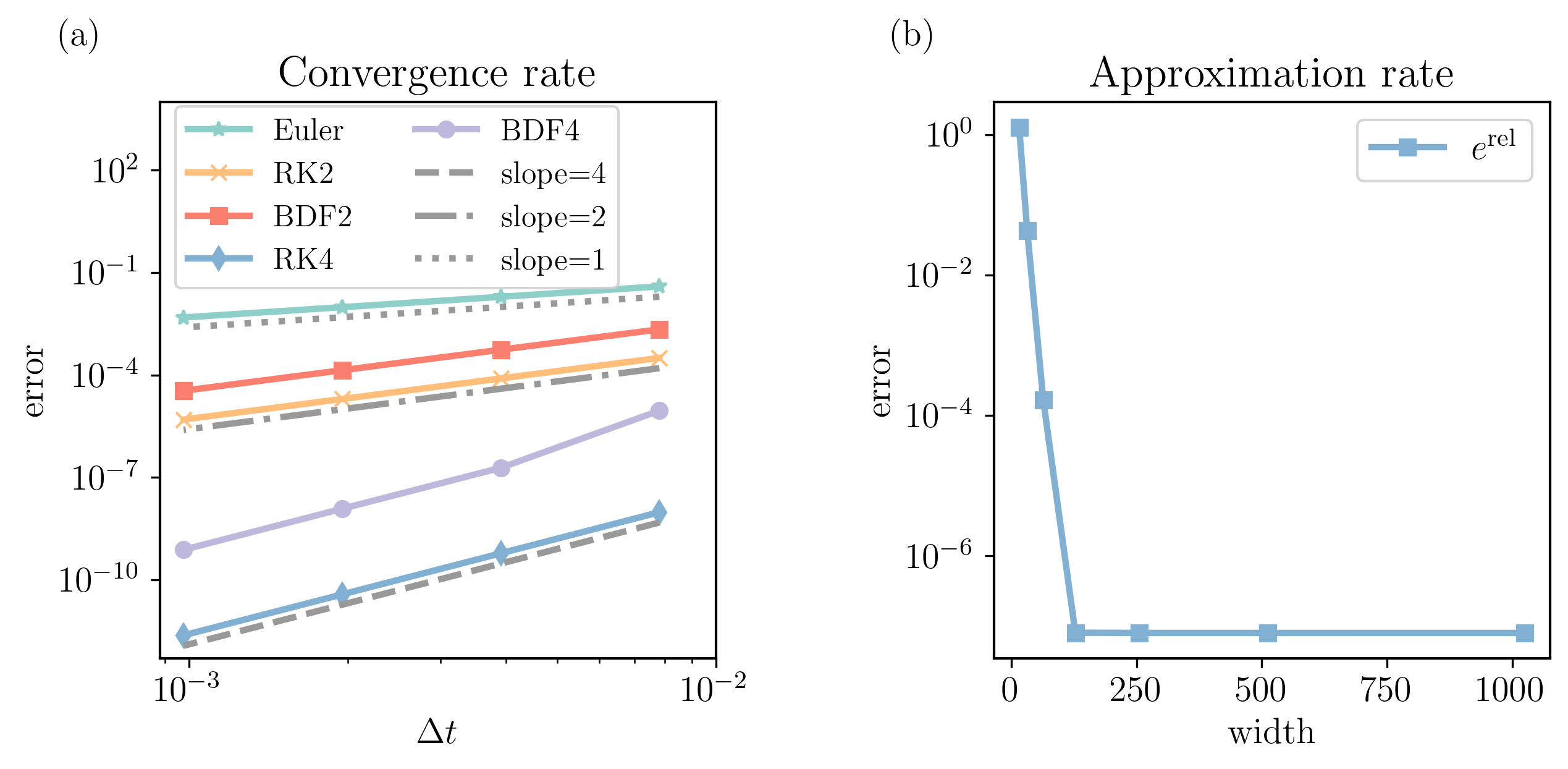}  
    \caption{Convergence rate tests for advection equation: (a) The slope of error curves for the five integrators can match their theoretical orders respectively, illustrating that the integrator within our framework remains convergence-guaranteed. (b) It can be observed that the error initially decreases at an exponential rate as the network width increases, as evidenced by the logarithmic scale on the y-axis. However, once the error reaches a certain level, it remains nearly constant, since it is then primarily determined by the accuracy of the time integrator.}
    \label{fig:converge}
\end{figure}

In this section, we use the advection equation as a representative example. As is well established in classic numerical analysis \cite{ErnstHairerSolving1993}, the Euler method achieves first-order accuracy, while both BDF2 and RK2 are second-order accurate, and BDF4 and RK4 attain fourth-order accuracy when spatial errors are neglected. Accordingly, we adopt a network of width $1000$ and examine the relative $L_2$ error at the final time step for various time step sizes. As shown in Fig.~\ref{fig:converge}(a), both the axes are both on a logarithmic scale, so that the slope of the error curve reflects the order of the accuracy. We can observe that all the five methods closely match the theoretical convergence orders. 

On the other hand, RNB-type networks are claimed to exhibit spectral accuracy \cite{DongLocal2021CMAME,JingrunChenBridging2022JML}. To assess this property, we evaluate the relative $L_2$ error of SDTM at the final time under varying network widths. In this experiment, the RK4 method is employed with $\Delta t = 0.001$, and the y-axis is plotted on a logarithmic scale. The results in Fig.~\ref{fig:converge}(b) show that the curve of error initially decreases linearly with increasing network width, indicating spectral accuracy within this regime. However, as the network width continues to grow, the error plateaus, suggesting that the spatial discretization error has become sufficiently small and the overall error is predominantly governed by the time integrator.

\subsection{Comparison of different integrators}
\begin{figure}[htpb]
    \centering
    \includegraphics[width=\textwidth]{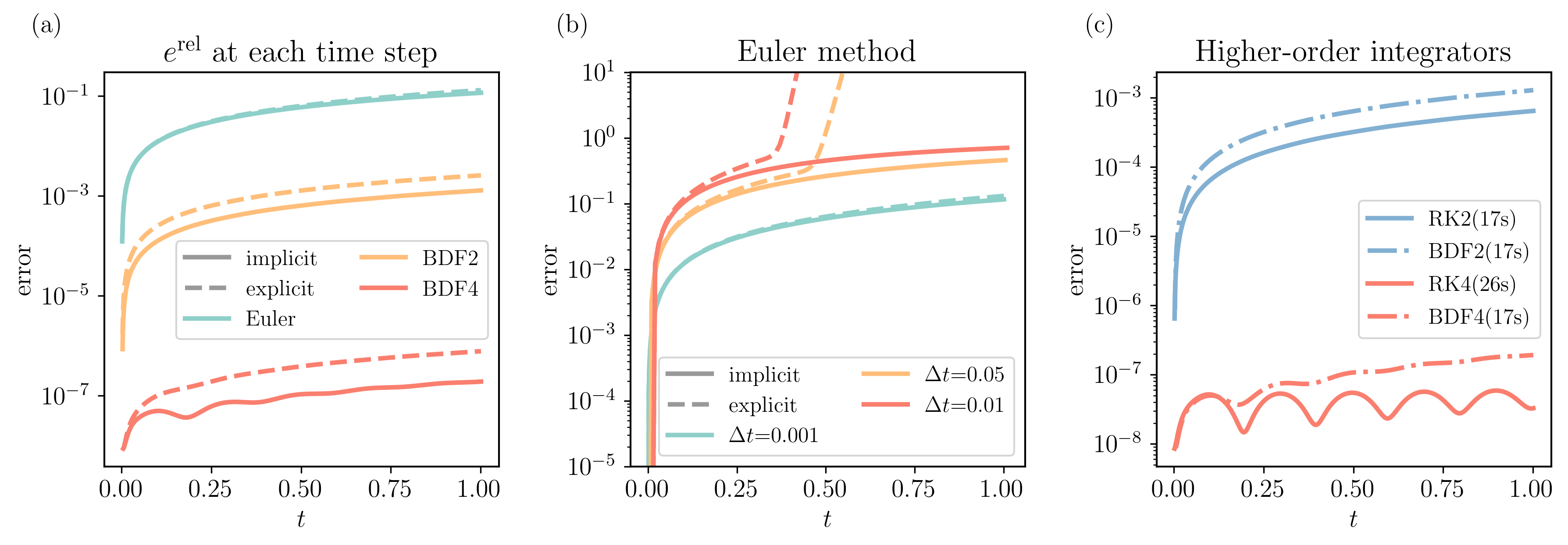}  
    \caption{Performance on different integrators for advection equation: (a) The performance of standard BDF method is better for their explicit variants; (b) The explicit Euler method is unstable when $\Delta t= 0.005$ and $0.01$, while the implicit Euler method is stable; (c) The RK type methods outperform the BDF type methods but are more time-consuming in higher order case.}
    \label{fig:integrator}
\end{figure}

In this section, we examine the differences among various time integrators, taking the advection equation as an illustrative example again. Our primary focus is on the comparison between explicit and implicit schemes. In this work, implicit methods are not employed for two main reasons: (1) when applied to nonlinear equations, implicit methods require solving nonlinear systems at each step, which is computationally demanding and makes it difficult to guarantee accuracy, often necessitating the use of semi-implicit schemes as a compromise; and (2) each iteration of an implicit method essentially involves solving an elliptic partial differential equation, which complicates the process of adaptively identifying a suitable initialization. Therefore, we often reformulate higher-order BDF methods in an explicit manner, i.e., Eq.~\eqref{eq:bdf2-ex} for BDF2. As shown in Fig.~\ref{fig:integrator}(a), the explicit BDF schemes generally exhibit slightly lower accuracy compared to their implicit counterparts, due to the incorporation of additional discretization errors. For the Euler method, the explicit and implicit variants perform nearly identically. In Fig.~\ref{fig:integrator}(b), we examine the performance of explicit and implicit Euler methods under different time step sizes. It can be observed that the explicit method becomes unstable when the time step is relatively large, whereas the implicit method maintains stability, consistent with classical results in numerical analysis. Finally, we compare the performance of different higher-order integrators, specifically the BDF and RK families. The BDF methods achieve higher-order accuracy by interpolating multiple historical solutions, whereas RK methods employ multi-stage slope evaluations within each integration step followed by weighted combinations. The results in Fig.~\ref{fig:integrator}(c) indicate that RK methods generally achieve higher accuracy than BDF methods. However, it is worth highlighting that while the computational cost of BDF2 is comparable to that of RK2, RK4 requires significantly more computation time, whereas BDF4 remains nearly as efficient as the second-order variant. The primary reason lies in the use of automatic differentiation: RK4 necessitates multiple derivative evaluations, which are more computationally expensive. This discrepancy in efficiency becomes even more pronounced for higher-order equations as shown in Tab.\ref{tab:cost}.

\subsection{Computation cost}
\label{sec:cost}
The computation cost of SDTM at each timestep is originating from two parts: generating the coefficient equations $A\theta \cong b$ for the least-squares problem and solving this problem. When the hidden-layers of the network remain invariant, all the time steps share the same $A$. Then we can precompute $A$ and its QR decomposition which helps the solve of least-squares problem. By this way, we only have to compute $b$ and solve $A\theta \cong b$ at each step, leading to a substantial improvement in SDTM`s efficiency. We report the execution time for different equations and integrators in Tab.~\ref{tab:cost}. Among all experiments, we can observe the following points.
\begin{enumerate}
    \item Precomputing the coefficient matrix can effectively enhance computational efficiency: the experiment labeled ``REINIT" represents that we reinitialize the network at each timestep. While we have to calculate the coefficient matrix repeatedly, this hampers the efficiency of computation.
    \item For higher order equation such as Burgers equation, while employing RK4 integrator, which requires differentiating the second-order equation four times, it takes more considerable computational cost and heavy GPU memory allocation in practice. Thus for higher order equations, we are used to choose higher order multi-step methods which only require a single differentiation.
    \item The computation cost is mainly dependent on number of collocation points $N$ and number of trainable parameters $M$. Since we can use non-uniform sampling method such as LHS, the required number of sampling points does not grow substantially as the dimensionality increases. Therefore, the computational efficiency in high-dimensional settings is comparable to that in low-dimensional cases.
\end{enumerate}
\begin{table}[h]
    \centering  
    \renewcommand{\arraystretch}{1.25}
    \begin{tabular}{c|c|c|c|c}
    \specialrule{1.2pt}{0pt}{0pt}
   \textbf{Equation}& \textbf{Setting(iterations)} & \textbf{Runtime(s)} &\textbf{Setting} &\textbf{Runtime(s)} \\
    \hline
      Advection & RK2 ($1\mathrm{e}4$)  & 79 &  RK4 & 130 \\
      \hline
       Burgers & BDF4 MsRNB($1\mathrm{e}4$) & 323 & RK4 & 6004 \\
               & BDF4 Adapt($1\mathrm{e}4$) & 341 & BDF4 REINIT & 10000+ \\ 
       \hline
      1D Allen--Cahn \eqref{eq:ac1d2} & BDF4($5\mathrm{e}3$) & 326 & MsRNB & 154 \\
       \hline 
      1D Allen--Cahn \eqref{eq:ac1d} & BDF4 Adapt($1\mathrm{e}4$) & 765 & MsRNB & 788 \\
       \hline
       2D Allen--Cahn & RK2($1\mathrm{e}3$) & 90 & Euler & 38\\
       & RK2($1\mathrm{e}4$) & 781 & Euler&349\\
       \hline
       Navier--Stokes & BDF4 ($1\mathrm{e}4$) & 1500 & BDF2 & 631\\
       \specialrule{1.2pt}{0pt}{0pt}
    \end{tabular}
    \caption{Time cost for different experiments. The iteration counts are omitted in the second column, since they are identical to those in the preceding column of the same row. The BDF4 method used here is an explicit version which is explained in Appendix.}
    \label{tab:cost}
\end{table}

\section{Concluding Remarks}
\label{sec:end}
In this work, we propose a semi-discrete in time framework towards time-dependent PDEs using random neural basis. We observe that temporal discretization substantially mitigates the approximation demands on the network relative to fully space–time coupled formulations. The temporal discretization is handled by classical numerical integrators, advancing the solution at discrete time steps, thereby inheriting the well-characterized numerical properties and interpretability of traditional methods. For the optimization at each timestep, the hidden-layer parameters of the network are fixed reducing each optimization to a linear least-squares problem, improving both computational efficiency and stability. For multi-scale PDEs, adaptive strategies are incorporated to adjust the network initialization or structure, enabling effective representation of features across scales. Numerical experiments not only demonstrate the effectiveness of the proposed method but also provide quantitative insights into the convergence rate of time integrators and the approximation rate of RNB. Overall, our framework effectively combines classical numerical analysis with the expressive power of neural networks, resulting in an efficient, interpretable, and accurate solver for time-dependent PDEs. 

In future work, we will conduct a systematic analysis of the stability of the SDTM framework. For example, in the case of conservation law equations, it is necessary to determine the corresponding CFL condition. However, since the present method does not involve spatial discretization, it is essential to further identify the conditions under which the network can deliver consistent performance. In terms of applications, we will extend the framework to higher-dimensional problems as well as those defined on irregular domains. Within this framework, the computational efficiency for such equations is not expected to degrade significantly; instead, the main challenge lies in designing neural networks that are well suited for these problems. As this direction aligns closely with the mainstream research trends in deep neural networks, we believe the proposed method holds strong potential for practical applications.

\section*{Acknowledgments}
This work is sponsored by the National Key R\&D Program of China Grant No. 2022YFA\-1008200 (T. L.) and Shanghai Institute for Mathematics and Interdisciplinary Sciences (T. L.).

\section*{Appendix}
\subsection*{A.Time Integrator}
\renewcommand{\theequation}{A.\arabic{equation}}
\setcounter{equation}{0} 
\label{sec:integrator}
We present some common time integrator schemes for evolutionary equations, which are thoroughly documented in standard numerical analysis literature \cite{ErnstHairerSolving1993}. Consider the following  semi-discrete formulation
\begin{equation}
    \vu_t = \mc F(\vu),
\end{equation}
where $\mc F$ is a differential operator for $\vx$. The choice of  time integrator falls naturally into two classes: one-step methods and multi-step methods. We present only the following subset of schemes, though most time integration schemes remain applicable within our framework. 

\noindent\textbf{One-step methods:}
The most common time integration scheme is explicit Euler method that reads
    \begin{equation}
        \vu^{n+1} = \vu^{n} + \Delta t\mc F(\vu^{n}).
    \end{equation}
This method is easy for implement but has only first order accuracy in time. For higher order methods, we employ Runge-Kutta methods, such as RK2 method
\begin{equation}
    \begin{aligned}
    \vk_1 &= \mc F(\vu^n),  \vk_2 = \mc F(\vu^n + \Delta t \vk_1),\\
    \vu^{n+1} &= \vu^{n} +  \frac{\Delta t}{2} (\vk_1 + \vk_2), 
    \end{aligned}
\end{equation}
and RK4 method
\begin{equation}
    \begin{aligned}
        \vk_1 &= \mc F(\vu^n), \quad \vk_2 = \mc F(\vu^n +  \frac{\Delta t}{2} \vk_1 ), \\ 
        \vk_3 &= \mc F(\vu^n + \frac{\Delta t}{2} \vk_2), \vk_4 = \mc F(\vu^n + \Delta t \vk_3) ,\\
        \vu^{n+1} &= \vu^n + \frac{\Delta t}{6}(\vk_1+ 2\vk_2+3\vk_3+\vk_4).
    \end{aligned}
\end{equation}
Despite the expansive parameter space of Runge-Kutta coefficient combinations, the present study exclusively employs the two time integration schemes specified above.

\noindent \textbf{Multi-step methods:}
In contrast to one-step schemes, multi-step methods necessitate a history of solution values to update $\vu^{n+1}$. This extended state dependency permits higher-order polynomial interpolation, achieving higher-order accuracy. We refer to backward differentiation formulas (BDF) methods, where a $k$-th order BDF method reads
\begin{equation}\label{eq:bdf}
    \sum_{j=1}^{k} \frac{1}{j}\nabla^j \vu^{n+1} = \Delta t \mc F(\vu^{n+1}), 
\end{equation}
here $\nabla^j$ is the backward difference quotient, i.e. $\nabla^1 \vu^{k+1} = \vu^{k+1} - \vu^{k}, \nabla^{j} \vu^{k+1} = \nabla^{j-1} \vu^{k+1} - \nabla^{j-1}\vu^{k}$ for $j>1$. Note that the right hand side of the Eq.~\eqref{eq:bdf} takes value at $t_{n+1}$, thus we get a class of implicit methods. When $k=1$, the BDF method is equivalent to explicit Euler method,
\begin{equation}
    \vu^{n+1} = \vu^{n} + \Delta t\mc F(\vu^{n+1}).
\end{equation}
When $k=2$ and $k=4$, we obtain BDF2 and BDF4 methods respectively
\begin{equation}
    \frac{3}{2} \vu^{n+1} - 2\vu^{n} + \frac{1}{2} \vu^{n-1} = \Delta t \mc F(\vu^{n+1}), 
\end{equation}

\begin{equation}
    \frac{25}{12} \vu^{n+1} - 4 \vu^{n} + 3\vu^{n-1} -\frac{4}{3}\vu^{n-2} + \frac{1}{4} \vu^{n-3} = \Delta t \mc{F} (\vu^{n+1}).
\end{equation}
However, for nonlinear problems, fully implicit schemes incur prohibitive computational cost due to per-timestep nonlinear solves. This motivates the adoption of implicit-explicit (IMEX) schemes that explicitly treat the nonlinear components while implicitly handle the linear operators. Ref.~\cite{HuangNew2024SJNA, HuangStability2023SJNAa} introduce a family of generalized BDF schemes with adapted IMEX implementations, which achieve enhanced numerical stability. Take $k=2$ for illustration, the origin method takes the backward differentiation operation at $t=t_{n+1}$. Furthermore, if we take this differentiation at $t=t_{n+\beta}$, $\beta>1$,  that we denote as BDF2-$\beta$ method, then the coefficients varies
\begin{equation}
    \frac{2\beta+1}{2} \vu^{n+1} - 2\beta\vu^{n} + \frac{2\beta-1}{2} \vu^{n-1} = \Delta t \mc F(\vu^{n+\beta}) + O(\Delta t^3), 
\end{equation}
then we can approximate the right-hand side with the history solutions explicitly 
\begin{equation}\label{eq:bdf2-ex}
    (\beta +1)\vu^{n} - \beta \vu^{n-1} = \vu^{n+\beta} + O(\Delta t^2),
\end{equation}
or implicitly
\begin{equation}\label{eq:bdf2-im}
    \beta \vu^{n+1} - (\beta -1)\vu^{n} = u^{n+\beta} + O(\Delta t^2).
\end{equation}
We denote this approach as explicit BDF$k$-$\beta$ (Ex-BDF$k$-$\beta$) method when adopting BDF$k$ integration at $t_{n+\beta}$ while interpolating the right-hand side with explicit formula \eqref{eq:bdf2-ex}. For the interpolations for $k>2$, one can easily derive the results using Taylor expansion and the detailed formulations are list in Ref.~\cite{HuangNew2024SJNA}.

\subsection*{B. Proof for Proposition 4.1}
\renewcommand{\theequation}{B.\arabic{equation}}
\setcounter{equation}{0} 
We briefly outline our proof. We use a complex analysis framework. We first compute the singular points of the basis functions $f_{\bm{k}}$. Then, using the Cauchy integral theorem, we find the decay rate of the Fourier coefficients and determine the constant $C$.
\begin{enumerate}[leftmargin=*]
\item Let $f_{\bm{k}} (\bm{x}) = \tanh(\bm{k} \cdot (\sin x_1, \dots, \sin x_n) )$ with $(x_1, \dots, x_n) \in \mathbb{R}^n$, we find that the singular points of $\tanh z = \frac{\sinh z }{\cosh z}$ locates at where $\cosh z = 0$. That is $z = i(\frac{\pi}{2} + m \pi)$. As a result, the singular points of $f_k$ satisfy the equation
\begin{equation}
\label{eq::5.12}
    k_1 \sin \xi_1 + \dots + k_n \xi_n  = i \left(\frac{\pi}{2} + m \pi \right),
\end{equation}
where $\xi_i = x_i + i y_i$. Using trigonometric identities, we have 
\begin{equation}
\label{eq::5.13}
    \sin (x+ iy) = \sin x \cosh y + i \cos x \sinh y.
\end{equation}
To find the singular point closest to the real plane, we consider $m=0$ and bring Eq.~\eqref{eq::5.13} into Eq.~\eqref{eq::5.12}.
\begin{equation}
    \sum_{i=1}^n k_i \left( \sin x_i \cosh y_i + i \cos x_i \sinh y_i \right) = i \frac{\pi}{2}.
\end{equation}
As a result, we get a sufficient condition to avoid singular points:
\begin{equation}
    \sum_{i=1}^n |k_i| |\sinh y_i| < \frac{\pi}{2}.
\end{equation}

\item We hope to choose a suitable integration path and obtain a consistent upper bound of $f_{\bm{k}}$ on this path that is independent of $\bm{k}$. We take 
\begin{equation}
    a_0 = \frac{1}{2} \operatorname{asinh} (\frac{\pi}{2 \| \bm{k}\|_1}),
\end{equation}
and get 
\begin{equation}
    \sum_{i=1}^n |k_i| \sinh a_0 = \| \bm{k}\|_1 \sinh a_0 \le  \| \bm{k}\|_1 \frac{\sinh  \operatorname{asinh} \frac{\pi}{2 \| \bm{k}\|_1}}{2} = \frac{\pi}{4}.
\end{equation}
Therefore, we obtain an integral path that avoids singularities. Moreover, since there is a uniform distance away from the starting point for any $\bm{k}$, we obtain an upper bound on the integral path for the function $f$ that is independent of $\bm{k}$:
\begin{equation}
    M_0:=\sup _{\operatorname{dist}(z, i \pi / 2) \geq \pi / 4}|\tanh z|<\infty .
\end{equation}

\item We define Fourier coefficient as
\begin{equation}
    \hat{f_{\bm{k}}}_{\bm{m}} = \frac{1}{(2 \pi)^n} \int_{[0,2\pi]^n} f_{\bm{k}} e^{-i \bm{m} \cdot \bm{x}} \operatorname{d} \bm{x}.
\end{equation}
Translate the $j$- th integral along the rectangular path upward ($m_j > 0$) or or downward ($m_j < 0$) to $\operatorname{sgn}(m_j) a_0$. By Cauchy's theorem and analyticity, we get
\begin{equation}
    \left|f_{\mathbf{m}}\right| \leq M_0 e^{-a_0 \sum_{j=1}^n\left|m_j\right|}
\end{equation}
To maintain $\hat{f_{\bm{k}}}_{\bm{m}} \le \varepsilon$, a sufficient condition is 
\begin{equation}
    \sum_{i=1}^n |m_i| \ge \frac{1}{a_0} \log \frac{M_0}{\varepsilon}.
\end{equation}
To sum up, the radius is 
\begin{equation}
    N_{\varepsilon} \approx \frac{1}{a_0} \log \frac{M_0}{\varepsilon} \sim \frac{2 \| \bm{k}\|_1}{\pi} \log \frac{M_0}{\varepsilon} .
\end{equation}

\end{enumerate}


\bibliography{ref}

\begin{thebibliography}{10}

\bibitem{AnagnostopoulosResidualbased2024CMAME}
Sokratis~J. Anagnostopoulos, Juan~Diego Toscano, Nikolaos Stergiopulos, and George~Em Karniadakis.
\newblock Residual-based attention in physics-informed neural networks.
\newblock {\em Computer Methods in Applied Mechanics and Engineering}, 421:116805, 2024.

\bibitem{BarronUniversal1993ITIT}
Andrew~R. Barron.
\newblock Universal approximation bounds for superpositions of a sigmoidal function.
\newblock {\em IEEE Transactions on Information Theory}, 39(3):930--945, 1993.

\bibitem{CalabroTime2023AMC}
Francesco Calabr{\`o}, Salvatore Cuomo, Daniela Di~Serafino, Giuseppe Izzo, and Eleonora Messina.
\newblock Time discretization in the solution of parabolic pdes with anns.
\newblock {\em Applied Mathematics and Computation}, 458:128230, 2023.

\bibitem{ChenDuality2023}
Hongrui Chen, Jihao Long, and Lei Wu.
\newblock A duality framework for generalization analysis of random feature models and two-layer neural networks, 2023.

\bibitem{ChenMicromacro2024}
Jingrun Chen, Zheng Ma, and Keke Wu.
\newblock A micro-macro decomposition-based asymptotic-preserving random feature method for multiscale radiative transfer equations, 2024.

\bibitem{ChenTENG2024P4ICML}
Zhuo Chen, Jacob McCarran, Esteban Vizcaino, Marin Solja{\v c}i{\'c}, and Di~Luo.
\newblock Teng: time-evolving natural gradient for solving pdes with deep neural nets toward machine precision.
\newblock In {\em Proceedings of the 41st International Conference on Machine Learning}, volume 235 of {\em ICML'24}, pages 7143--7162, Vienna, Austria, 2024. JMLR.org.

\bibitem{CiarletFinite2002}
Philippe~G. Ciarlet.
\newblock {\em The finite element method for elliptic problems}.
\newblock {Society for Industrial and Applied Mathematics}, 2002.

\bibitem{DonatellaDynamics2023PRA}
Kaelan Donatella, Zakari Denis, Alexandre Le~Boit{\'e}, and Cristiano Ciuti.
\newblock Dynamics with autoregressive neural quantum states: Application to critical quench dynamics.
\newblock {\em Physical Review A}, 108(2):22210, 2023.

\bibitem{DongLocal2021CMAME}
Suchuan Dong and Zongwei Li.
\newblock Local extreme learning machines and domain decomposition for solving linear and nonlinear partial differential equations.
\newblock {\em Computer Methods in Applied Mechanics and Engineering}, 387:114129, 2021.

\bibitem{DongMethod2021JCP}
Suchuan Dong and Naxian Ni.
\newblock A method for representing periodic functions and enforcing exactly periodic boundary conditions with deep neural networks.
\newblock {\em Journal of Computational Physics}, 435:110242, 2021.

\bibitem{DuffyFinite2006}
Daniel~J. Duffy.
\newblock {\em Finite difference methods in financial engineering: a partial differential equation approach}.
\newblock Wiley finance series. John Wiley, Chichester, England ; Hoboken, NJ, 2006.

\bibitem{ErnstHairerSolving1993}
{Ernst Hairer}, {Gerhard Wanner Wanner}, and {Syvert P. N{\o}rsett}.
\newblock {\em Solving ordinary differential equations I}.
\newblock Springer Series in Computational Mathematics. Springer Berlin Heidelberg, Berlin, Heidelberg, 1993.

\bibitem{FanMultiscale2019MMS}
Yuwei Fan, Lin Lin, Lexing Ying, and Leonardo {Zepeda-N{\'u}\~{n}ez}.
\newblock A multiscale neural network based on hierarchical matrices.
\newblock {\em Multiscale Modeling and Simulation}, 17(4):1189--1213, 2019.

\bibitem{HuangStability2023SJNAa}
Fukeng Huang and Jie Shen.
\newblock Stability and error analysis of a second-order consistent splitting scheme for the navier--stokes equations.
\newblock {\em SIAM Journal on Numerical Analysis}, 61(5):2408--2433, 2023.

\bibitem{HuangNew2024SJNA}
Fukeng Huang and Jie Shen.
\newblock On a new class of bdf and imex schemes for parabolic type equations.
\newblock {\em SIAM Journal on Numerical Analysis}, 62(4):1609--1637, 2024.

\bibitem{HuangFrequencyadaptive2025CMAME}
Jizu Huang, Rukang You, and Tao Zhou.
\newblock Frequency-adaptive multi-scale deep neural networks.
\newblock {\em Computer Methods in Applied Mechanics and Engineering}, 437:117751, 2025.

\bibitem{JingrunChenBridging2022JML}
Jingrun~Chen Jingrun~Chen, Xurong~Chi Xurong~Chi, Weinan~E Weinan~E, and Zhouwang~Yang Zhouwang~Yang.
\newblock Bridging traditional and machine learning-based algorithms for solving pdes: the random feature method.
\newblock {\em Journal of Machine Learning}, 1(3):268--298, 2022.

\bibitem{KamontNumerical1999HFDIaA}
Zdzislaw Kamont.
\newblock Numerical method of lines.
\newblock In Zdzislaw Kamont, editor, {\em Hyperbolic Functional Differential Inequalities and Applications}, pages 181--204. Springer Netherlands, Dordrecht, 1999.

\bibitem{LiMultiscale2020C}
Xi-An Li, Zhi-Qin~John Xu, and Lei Zhang.
\newblock A multi-scale dnn algorithm for nonlinear elliptic equations with multiple scales.
\newblock {\em Communications in Computational Physics}, 28(5):1886--1906, 2020.

\bibitem{LiuMultiscale2020CCP}
Ziqi Liu, Wei Cai, and Zhi-Qin~John Xu.
\newblock Multi-scale deep neural network (mscalednn) for solving poisson-boltzmann equation in complex domains.
\newblock {\em Communications in Computational Physics}, 28(5):1970--2001, 2020.

\bibitem{LuoAutoregressive2022PRL}
Di~Luo, Zhuo Chen, Juan Carrasquilla, and Bryan~K Clark.
\newblock Autoregressive neural network for simulating open quantum systems via a probabilistic formulation.
\newblock {\em Physical Review Letters}, 2022.

\bibitem{MatteyNovel2022CMAME}
Revanth Mattey and Susanta Ghosh.
\newblock A novel sequential method to train physics informed neural networks for allen cahn and cahn hilliard equations.
\newblock {\em Computer Methods in Applied Mechanics and Engineering}, 390:114474, 2022.

\bibitem{PoochinapanNumerical2022AMC}
Kanyuta Poochinapan and Ben Wongsaijai.
\newblock Numerical analysis for solving allen-cahn equation~in 1d and 2d based on higher-order compact structure-preserving difference scheme.
\newblock {\em Applied Mathematics and Computation}, 434:127374, 2022.

\bibitem{RaissiPhysicsinformed2019JCP}
M.~Raissi, P.~Perdikaris, and G.E. Karniadakis.
\newblock Physics-informed neural networks: a deep learning framework for solving forward and inverse problems involving nonlinear partial differential equations.
\newblock {\em Journal of Computational Physics}, 378:686--707, 2019.

\bibitem{SchatzQuantum2002}
George~C. Schatz and Mark~A. Ratner.
\newblock {\em Quantum mechanics in chemistry}.
\newblock Dover Publ, Mineola, NY, dover ed., corr., unabridged republ edition, 2002.

\bibitem{ShangRandomized2024JEM}
Yong Shang and Fei Wang.
\newblock Randomized neural networks with petrov--galerkin methods for solving linear elasticity and navier--stokes equations.
\newblock 150(4):04024010.

\bibitem{SommerfeldIntroduction1949PDEiP}
Arnold Sommerfeld.
\newblock Introduction to partial differential equations.
\newblock In {\em Partial Differential Equations in Physics}, pages 32--62. Elsevier, 1949.

\bibitem{StoerIntroduction2002}
J.~Stoer and R.~Bulirsch.
\newblock {\em Introduction to numerical analysis}, volume~12 of {\em Texts in Applied Mathematics}.
\newblock Springer New York, New York, NY, 2002.

\bibitem{SunLocal2024JCAM}
Jingbo Sun, Suchuan Dong, and Fei Wang.
\newblock Local randomized neural networks with discontinuous galerkin methods for partial differential equations.
\newblock 445:115830.

\bibitem{SunSurrogate2020CMAME}
Luning Sun, Han Gao, Shaowu Pan, and Jian-Xun Wang.
\newblock Surrogate modeling for fluid flows based on physics-constrained deep learning without simulation data.
\newblock {\em Computer Methods in Applied Mechanics and Engineering}, 361:112732, 2020.

\bibitem{WangPirateNets2024JMLR}
Sifan Wang, Bowen Li, Yuhan Chen, and Paris Perdikaris.
\newblock Piratenets: Physics-informed deep learning with residual adaptive networks.
\newblock {\em Journal of Machine Learning Research}, 25(402):1--51, 2024.

\bibitem{WangUnderstanding2021SJSC}
Sifan Wang, Yujun Teng, and Paris Perdikaris.
\newblock Understanding and mitigating gradient flow pathologies in physics-informed neural networks.
\newblock {\em SIAM Journal on Scientific Computing}, 43(5):A3055--A3081, 2021.

\bibitem{WangWhen2022JCP}
Sifan Wang, Xinling Yu, and Paris Perdikaris.
\newblock When and why pinns fail to train: a neural tangent kernel perspective.
\newblock {\em Journal of Computational Physics}, 449:110768, 2022.

\bibitem{WangExtreme2024CMAME}
Yiran Wang and Suchuan Dong.
\newblock An extreme learning machine-based method for computational pdes in higher dimensions.
\newblock {\em Computer Methods in Applied Mechanics and Engineering}, 418:116578, 2024.

\bibitem{ZhangTransferable2024JSC}
Zezhong Zhang, Feng Bao, Lili Ju, and Guannan Zhang.
\newblock Transferable neural networks for partial differential equations.
\newblock {\em Journal of Scientific Computing}, 99(1):2, 2024.

\end{thebibliography}
\bibliographystyle{plain}

\end{document}